\documentclass[11pt]{article}
\usepackage{amsfonts,amsmath,amsxtra}
\usepackage{latexsym}
\usepackage{amssymb}

\def\hybrid{\topmargin 0pt      \oddsidemargin 0pt
        \headheight 0pt \headsep 0pt
        \textwidth 16.5cm
        \textheight 23cm
        \marginparwidth 0.0in
        \parskip 5pt plus 1pt   \jot = 1.5ex}
\catcode`\@=11
\def\marginnote#1{}
\newcount\hour
\newcount\minute
\newtoks\amorpm
\hour=\time\divide\hour by60 \minute=\time{\multiply\hour by60
\global\advance\minute by-\hour}
\edef\standardtime{{\ifnum\hour<12 \global\amorpm={am}%
        \else\global\amorpm={pm}\advance\hour by-12 \fi
        \ifnum\hour=0 \hour=12 \fi
      \number\hour:\ifnum\minute<10 0\fi\number\minute\the\amorpm}}
\edef\militarytime{\number\hour:\ifnum\minute<10 0\fi\number\minute}

\def\draftlabel#1{{\@bsphack\if@filesw {\let\thepage\relax
   \xdef\@gtempa{\write\@auxout{\string
      \newlabel{#1}{{\@currentlabel}{\thepage}}}}}\@gtempa
   \if@nobreak \ifvmode\nobreak\fi\fi\fi\@esphack}
        \gdef\@eqnlabel{#1}}
\def\@eqnlabel{}
\def\@vacuum{}
\def\draftmarginnote#1{\marginpar{\raggedright\scriptsize\tt#1}}

\def\draft{\oddsidemargin -0.1truein
        \def\@oddfoot{\sl preliminary draft \hfil
        \rm\thepage\hfil\sl\today\quad\militarytime}
        \let\@evenfoot\@oddfoot \overfullrule 3pt
        \let\label=\draftlabel
        \let\marginnote=\draftmarginnote
\def\@eqnnum{{\rm (\theequation)}
\rlap{\kern\marginparsep\tt\@eqnlabel}%
\global\let\@eqnlabel\@vacuum}  }

\newfont{\Bbbb}{msbm7 scaled 1\@ptsize00}
\newcommand{\zs}{\raise-1pt\hbox{$\mbox{\Bbbb Z}$}}

scaled 1\@ptsize00

\font\sevenmsa=msam6 
\newfam\msafam
\textfont\msafam=\sevenmsa
\def\hexnumber@#1{\ifnum#1<10 \number#1\else
\ifnum#1=10 A\else\ifnum#1=11 B\else\ifnum#1=12 C\else \ifnum#1=13
D\else\ifnum#1=14 E\else\ifnum#1=15 F\fi\fi\fi\fi\fi\fi\fi}
\def\msa@{\hexnumber@\msafam}
\def\llcorner{\delimiter"4\msa@78\msa@78 }
\def\lrcorner{\delimiter"5\msa@79\msa@79 }
\mathchardef\blacktriangleright="3\msa@49
\mathchardef\blacktriangleleft="3\msa@4A \font\tenmsb=msbm10 scaled
1\@ptsize00
\newfam\msbfam
\textfont\msbfam=\tenmsb \scriptfont\msbfam=\tenmsb

\newdimen\Squaresize \Squaresize=14pt
\newdimen\Thickness \Thickness=0.5pt

\def\Square#1{\hbox{\vrule width \Thickness
   \vbox to \Squaresize{\hrule height \Thickness\vss
      \hbox to \Squaresize{\hss#1\hss}
   \vss\hrule height\Thickness}
\unskip\vrule width \Thickness} \kern-\Thickness}

\def\Vsquare#1{\vbox{\Square{$#1$}}\kern-\Thickness}

\def\numberbysection{\@addtoreset{equation}{section}
        \def\theequation{\thesection.\arabic{equation}}}
\numberbysection

\renewcommand{\theequation}{\thesection.\arabic{equation}}
\def\titlepage{\@restonecolfalse\if@twocolumn\@restonecoltrue\onecolumn
     \else \newpage \fi \thispagestyle{empty}\c@page\z@
        \def\thefootnote{\fnsymbol{footnote}} }

\def\endtitlepage{\if@restonecol\twocolumn \else  \fi
        \def\thefootnote{\arabic{footnote}}
        \setcounter{footnote}{0}}  
\relax

\hybrid
\makeatletter
\newdimen\normalarrayskip            
\newdimen\minarrayskip               
\normalarrayskip\baselineskip \minarrayskip\jot
\newif\ifold             \oldtrue            \def\new{\oldfalse}
\def\arraymode{\ifold\relax\else\displaystyle\fi}
\def\eqnumphantom{\phantom{(\theequation)}} 
\def\@arrayskip{\ifold\baselineskip\z@\lineskip\z@
     \else
     \baselineskip\minarrayskip\lineskip1\baselineskip\fi}

\def\@arrayclassz{\ifcase \@lastchclass \@acolampacol \or
\@ampacol \or \or \or \@addamp \or
   \@acolampacol \or \@firstampfalse \@acol \fi
\edef\@preamble{\@preamble
  \ifcase \@chnum
     \hfil$\relax\arraymode\@sharp$\hfil
     \or $\relax\arraymode\@sharp$\hfil
     \or \hfil$\relax\arraymode\@sharp$\fi}}

\def\@array[#1]#2{\setbox\@arstrutbox=\hbox{\vrule
     height\arraystretch \ht\strutbox
     depth\arraystretch \dp\strutbox
width\z@}\@mkpream{#2}\edef\@preamble{\halign \noexpand\@halignto
\bgroup \tabskip\z@ \@arstrut \@preamble \tabskip\z@ \cr}%
\let\@startpbox\@@startpbox \let\@endpbox\@@endpbox
  \if #1t\vtop \else \if#1b\vbox \else \vcenter \fi\fi
  \bgroup \let\par\relax
  \let\@sharp##\let\protect\relax
  \@arrayskip\@preamble}
\def\eqnarray{\stepcounter{equation}%
              \let\@currentlabel=\theequation
              \global\@eqnswtrue
              \global\@eqcnt\z@
              \tabskip\@centering              
              \let\\=\@eqncr
              $$%
            \halign to \displaywidth  \bgroup
             \eqnumphantom \@eqnsel
      \hskip\@centering                               
    $\displaystyle  \tabskip\z@ {##}$%
    &\global\@eqcnt\@ne \hskip 2\arraycolsep
         $ \displaystyle  \arraymode{##}$\hfil
    &\global\@eqcnt\tw@ \hskip 2\arraycolsep
         $\displaystyle\tabskip\z@{##}$\hfil
         \tabskip\@centering
    &{##}\tabskip\z@\cr}
\makeatother
\newcommand{\RR}{{\mathbb{R}}}

\newcommand{\ZZ}{{\mathbb{Z}}}

\newcommand{\C}{\mathbb{C}}

\def\IC{\mathbb{C}}

\def\IP{\mathbb{P}}
\def\IQ{\mathbb{Q}}

\def\IZ{\mathbb{Z}}

\def\CD {\mathcal{D}}
\def\CE {\mathcal{E}}

\def\CH {\mathcal{H}}

\def\CK {\mathcal{K}}
\def\CL {\mathcal{L}}
\def\CM {\mathcal{M}}
\def\CN {\mathcal{N}}
\def\CO {\mathcal{O}}
\def\CP {\mathcal{P}}
\def\CQ {\mathcal{Q}}
\def\CR {\mathcal{R}}

\def\CT {\mathcal{T}}

\def\CV {\mathcal{V}}

\def\QM{{\CQ\CM}}

\def\ch{{\cal H}}

\def\la{\lambda}

\def\pr {\partial}

\def\ch{{\rm ch}}

\def\Lie{{\rm Lie}}

\def\pt{{\rm pt}}

\def\Sym{{\mathop{\rm Sym}}}
\def\Tr{{\rm Tr}}

\newtheorem{te}{Theorem}[section]

\newtheorem{prop}{Proposition}[section]           
\newtheorem{cor}{Corollary}[section]
\newtheorem{lem}{Lemma}[section]

\newtheorem{rem}{Remark}[section]

\newtheorem{prob}{Porblem}[section]
\newcommand\bqa{\begin{eqnarray}}
\newcommand\eqa{\end{eqnarray}}
\def\be{\begin{eqnarray}\new\begin{array}{cc}}
\def\ee{\end{array}\end{eqnarray}}

\def\beq{\begin{equation}}
\def\eeq{\end{equation}}
\def\bse{\begin{subequations}}                
\def\ese{\end{subequations}}
\def\bp{\begin{pmatrix}}
\def\ep{\end{pmatrix}}

\def\i{\imath}

\def\stack#1#2{\raise0.7pt\hbox{$\mathrel{\mathop{#2}\limits^{#1}}$}}
\def\tr{\triangleright}
\def\tl{\triangleleft}
\def\sem{\mathsurround=0pt \raise1pt
\hbox{$\scriptscriptstyle>\!\!$}\:\!\!\tl}
\def\mes{\mathsurround=0pt \tr\!\:\!\raise0.8pt
\hbox{$\scriptscriptstyle\!\!<$}\,}
\def\]{\mathsurround=0pt ]\raise-2pt\hbox{$_\ast$}}

\def\<{\langle}
\def\>{\rangle}

\def\CQ{{\cal Q}}

\def\ch{{\cal H}}

\def\CO{{\cal O}}

\def\CH{\mathcal{H}}

\def\we{\raise-1pt\hbox{$\,\stackrel{\wedge}{,}\,$}}
\def\tr{{\rm tr}\,}
\def\Tr{{\rm Tr}\,}
\def\pr {\partial}

\newcounter{pac}[section]

0

0

\title{\bf On $q$-deformed $\mathfrak{gl}_{\ell+1}$-Whittaker function II}
\begin{document}
\author{Anton Gerasimov, Dimitri Lebedev, and Sergey Oblezin}
\date{}

\maketitle

\renewcommand{\abstractname}{}

\begin{abstract}
\noindent {\bf Abstract}. A representation of a specialization of
a $q$-deformed class one lattice $\mathfrak{gl}_{\ell+1}$-Whittaker
function in terms of
cohomology groups of line bundles on the space $\QM_d(\IP^{\ell})$
of quasi-maps $\IP^1\to \IP^{\ell}$ of degree $d$    is proposed. For
$\ell=1$,  this provides an interpretation of non-specialized
$q$-deformed $\mathfrak{gl}_{2}$-Whittaker function in terms of
$\QM_d(\IP^1)$. In particular the ($q$-version of) Mellin-Barnes
representation of  $\mathfrak{gl}_2$-Whittaker function is
realized as a  semi-infinite period map. The explicit form of the
period map  manifests an important role of $q$-version of
$\Gamma$-function as a substitute of topological genus
in semi-infinite geometry. 
A relation with  Givental-Lee universal solution ($J$-function) of $q$-deformed
$\mathfrak{gl}_2$-Toda chain is also discussed.

\end{abstract}

\section*{Introduction}

In the first part \cite{GLO1} of this series of papers we have proposed an
explicit representation of a  $q$-deformed class one lattice 
$\mathfrak{gl}_{\ell+1}$-Whittaker function defined as a common
eigenfunction of a complete set of commuting quantum Hamiltonians of
$q$-deformed $\mathfrak{gl}_{\ell+1}$-Toda chain. Here ``class one''
means that Whittaker function is non-zero only in the dominant
domain. On $q$-deformed Toda chains  see e.g.
\cite{Et}. The case $\ell=1$ was discussed previously in
\cite{GLO3} (for related results in this direction see \cite{KLS},
\cite{GiL}, \cite{GKL1}, \cite{BF}, \cite{FFJMM}).
A special feature  of the proposed representation is that
 $q$-deformed class one $\mathfrak{gl}_{\ell+1}$-Whittaker
function $\Psi^{\mathfrak{gl}_{\ell+1}}_{\underline{z}}(\underline{p})$ with
$\underline{z}=(z_1,\ldots,z_{\ell+1})$ and
$\underline{p}=(p_1,\ldots,p_{\ell+1})\in\IZ^{\ell+1}$, is given by  a
character of a  $\IC^*\times GL_{\ell+1}(\IC)$-module $\CV_{\underline{p}}$.
The expression in terms of a character can be considered as a $q$-version of
Shintani-Casselman-Shalika representation of class one $p$-adic Whittaker
functions \cite{Sh}, \cite{CS}. Indeed  our
representation of $q$-deformed $\mathfrak{gl}_{\ell+1}$-Whittaker
function reduces, in a certain limit, to the Shintani-Casselman-Shalika representation of
$p$-adic Whittaker function. Note that the representation
$q$-deformed Whittaker function  as a
character is a $q$-analog of  the Givental integral representation 
 \cite{Gi2}, \cite{GKLO} of the classical
$\mathfrak{gl}_{\ell+1}$-Whittaker function.

The main objective of this paper is a better  understanding of the
representation  of $q$-deformed $\mathfrak{gl}_{\ell+1}$-Whittaker
function as a  character.  Below we  will consider a specialization 
of the $q$-deformed Whittaker function given by the 
trace over $\IC^*\times GL_{\ell+1}(\IC)$-module $\CV_{n,k}$ (in the
case $\ell=1$ there is actually no specialization).
Our main result is  presented  in Theorem \ref{mainTh21}.
We provide a description of $\IC^*\times
GL_{\ell+1}(\IC)$-module $\CV_{n,k}$ as a zero degree cohomology group
of a line bundle on an algebraic version 
 $\CL\IP^{\ell}_+$ of 
a semi-infinite cycle $\widetilde{L\IP^\ell}_+$ in a universal covering
$\widetilde{L\IP^\ell}$ of the space of loops in $\IP^\ell$. We
define $\CL\IP^\ell_+$ as an appropriate  limit $d\to
\infty$ of the space $\QM_d(\IP^\ell)$ 
of degree $d$ quasi-maps of $\IP^1$ to
$\IP^\ell$  \cite{Gi1}, \cite{CJS}.
In particular for $\ell=1$ this provides a description of
a $q$-deformed $\mathfrak{gl}_2$-Whittaker function  in terms
of cohomology of line bundles over $\CL\IP^1_+$.
A universal solution of the  $q$-deformed
$\mathfrak{gl}_{\ell+1}$-Toda chain \cite{GiL} was given  in terms
of cohomology groups of line bundles over $\QM_d(X), X=G/B$ for
finite $d$. We  demonstrate how our interpretation of the  $q$-deformed
$\mathfrak{gl}_{\ell+1}$-Whittaker function is reconciled with the results
of \cite{GiL}.

Using Theorems \ref{mainTh21}, we
interpret a $q$-version of the Mellin-Barnes integral
representation of the specialized $q$-Whittaker function
as a semi-infinite analog of the Riemann-Roch-Hirzebruch theorem.
The corresponding Todd class is expressed in terms  
of a $q$-version of $\Gamma$-function.  
Analogously, the classical $\Gamma$-function appears in a description of  
the fundamental class of semi-infinite homology theory
and enters the Mellin-Barnes integral representation of 
the classical Whittaker function. We briefly consider an analog of the
elliptic genus arising in the $\IC^*$-localization on 
 $\CL\IP^{\ell}_+$. We demonstrate that proliferation
of fixed points of $\IC^*$-action obstructs identification of the result 
as a topological genus of an extraordinary  cohomology theory. 
Note also that  the ($q$-version of) $\Gamma$-function which appears in our
calculations of a semi-infinite analog of the Todd class was 
considered as a candidate for a topological genus  
by Kontsevich \cite{K} (see also \cite{Li},\cite{Ho}).

Let us stress that the $\IC^*\times GL_{\ell+1}(\IC)$-module $\CV_{n,k}$
arising in the description of $q$-deformed
$\mathfrak{gl}_{\ell+1}$-Whittaker function  is not
irreducible. It would be natural to look for an interpretation of
$\CV_{n,k}$ as an irreducible module of  a quantum affine Lie group.
 A relation of the geometry of semi-infinite flags
to representation theory of affine Lie algebras
was proposed  in \cite{FF}.
The semi-infinite flag space is defined as
 $X^{\frac{\infty}{2}}=G(\CK)/H(\CO)N(\CK)$ where $\CK=\IC((t))$,
$\CO=\IC[[t]]$, $B=HN$ is a Borel subgroup of $G$,  $N$ is its
unipotent radical and $H$ is the associated Cartan subgroup.
The semi-infinite flag spaces are not easy to deal with. 
An interesting approach to the semi-infinite geometry 
was proposed by Drinfeld. He introduced a space of quasi-maps
$\CQ\CM_{\underline{d}}(\IP^1,G/B)$ that 
should be considered as a finite-dimensional substitute of the
semi-infinite flag space $X^{\frac{\infty}{2}}$ (see e.g. \cite{FM},
\cite{FFM}, \cite{Bra}). Thus, taking into
account constructions proposed in this paper  one can expect that
($q$-deformed) $\mathfrak{gl}_{\ell+1}$-Whittaker functions
(encoding Gromov-Witten invariants and their $K$-theory generalizations)
can be expressed in terms of  representation theory of affine Lie
algebras (see \cite{GiL} for a related conjecture and 
\cite{FFJMM} for a recent progress in this direction).
The paper  \cite{GLO2}  deals with a relation 
of our results with the representation theory of (quantum) affine Lie groups.

The paper is organized as follows. In Section 1, 
explicit solutions of $q$-deformed $\mathfrak{gl}_{\ell+1}$-Toda chain
($q$-versions of Whittaker functions) are  recalled.
In Section 2,  we derive integral expressions for the 
counting of holomorphic sections of line bundles
in the space of quasi-maps.  
In Section 3 we derive  a representation of specialized
$q$-Whittaker functions in
terms of cohomology of holomorphic line bundles on the space of quasi-maps of
$\IP^1$ to $\IP^\ell$. We propose an interpretation of $q$-Whittaker
functions as semi-infinite periods.  
In Section 4 the analogous interpretation of the classical Whittaker
functions is discussed. In Section 5, we clarify the connection of our
interpretation of the  $q$-deformed $\mathfrak{gl}_{\ell+1}$-Whittaker function
with the results of \cite{GiL}.
Finally,  in Section 6 we consider an analog of elliptic genus
arising after $\IC^*$-localization on $\CL\IP^{\ell}_+$ 
 and its possible relation with extraordinary cohomology theories.

{\em Acknowledgments}: We are grateful to A.~Rosly for useful
discussions. The research of AG was  partly supported by
SFI Research Frontier Programme and Marie Curie RTN Forces Universe
from EU. The research of SO is partially supported by  RF President
Grant MK-134.2007.1.  The research was also partially supported by
Grant RFBR-08-01-00931-a.

\section{ $q$-deformed  $\mathfrak{gl}_{\ell+1}$-Whittaker
function}

In this section we recall a  construction \cite{GLO1} of the $q$-deformed
$\mathfrak{gl}_{\ell+1}$-Whittaker function
$\Psi^{\mathfrak{gl}_{\ell+1}}_{\underline{z}}
(\underline{p}_{\ell+1})$ defined on the lattice
$\underline{p}_{\ell+1}=(p_{\ell+1,1},\ldots, p_{\ell+1,\ell+1})
\in\ZZ^{\ell+1}$. We will consider only class one Whittaker functions
satisfying the condition
$$
\Psi^{\mathfrak{gl}_{\ell+1}}_{\underline{z}}
(\underline{p}_{\ell+1})=0
$$
outside dominant domain $p_{\ell+1,1}\geq\ldots\geq
p_{\ell+1,\ell+1}$.

The $q$-deformed $\mathfrak{gl}_{\ell+1}$-Whittaker functions are
common eigenfunctions of $q$-deformed $\mathfrak{gl}_{\ell+1}$-Toda
chain Hamiltonians:  \be\label{comm}
\CH_r^{\mathfrak{gl}_{\ell+1}}(\underline{p}_{\ell+1})\,=\,\sum_{I_r}\,\bigl(
\widetilde{X}_{i_1}^{1-\delta_{i_2-i_1,\,1}}\cdot\ldots\cdot
\widetilde{X}_{i_{r-1}}^{1-\delta_{i_r-i_{r-1},\,1}}\cdot
\widetilde{X}_{i_r}^{1-\delta_{i_{r+1}-i_r,\,1}}\bigr)
T_{i_1}\cdot\ldots\cdot T_{i_r},\ee where 
the sum is over ordered subsets $I_r=\{i_1<i_2<\ldots <i_r\}\subset
\{1,2,\ldots ,\ell+1\}$ and we assume
$i_{r+1}=\ell+2$. In \eqref{comm}  we use  the following notations
$$
T_if(\underline{p}_{\ell+1})=f(\underline{\widetilde{p}}_{\ell+1}),
\hspace{1.5cm}\widetilde{p}_{\ell+1,k}=p_{\ell+1,k}+\delta_{k,i},
\hspace{2cm}i,k=1,\ldots,\ell+1,
$$
$$
\widetilde{X}_i=1-q^{p_{\ell+1,i}-p_{\ell+1,i+1}+1},\hspace{0.5cm}
i=1,\ldots,\ell,
$$
and  $\widetilde{X}_{\ell+1}=1$. We  assume $q\in \IC^*,|q|<1$.
For example, the   first nontrivial Hamiltonian 
has the following form:
\be\label{FirstHamiltonian}
 {\CH}_1^{\mathfrak{gl}_{\ell+1}}(\underline{p}_{\ell+1})\,=\,
\sum\limits_{i=1}^{\ell}(1-q^{p_{\ell+1,i}-p_{\ell+1,i+1}+1})T_i\,+\,
T_{\ell+1}. \ee
 The main result of \cite{GLO1} is a construction of 
common eigenfunctions of quantum Hamiltonians \eqref{comm}:
 \be\label{eiglat}
\ch_r^{\mathfrak{gl}_{\ell+1}}(\underline{p}_{\ell+1})
\Psi^{\mathfrak{gl}_{\ell+1}}_{z_1,\ldots ,z_{\ell+1}}
(\underline{p}_{\ell+1})\,=\,(\sum_{ I_r}\prod\limits_{i\in I_r} z_i
)\,\,\Psi^{\mathfrak{gl}_{\ell+1}}_{z_1,\ldots ,z_{\ell+1}}
(\underline{p}_{\ell+1}).\ee Denote by
$\CP^{(\ell+1)}\subset\ZZ^{\ell(\ell+1)/2}$ a  subset of integers 
$p_{n,i}$, $n=1,\ldots,\ell+1$, $i=1,\ldots,n$ satisfying the
Gelfand-Zetlin conditions $p_{k+1,i}\geq p_{k,i}\geq p_{k+1,i+1}$
for $k=1,\ldots,\ell$.  In the following we use the standard
notation $(n)_q!=(1-q)...(1-q^n)$.

\begin{te}\label{mainth1}
Let
$\Psi^{\mathfrak{gl}_{\ell+1}}_{z_1,\ldots,z_{\ell+1}}(\underline{p}_{\ell+1})$
be a function given in the dominant  domain $p_{\ell+1,1}\geq\ldots\geq
p_{\ell+1,\ell+1}$ by
\be\label{main}
 \Psi^{\mathfrak{gl}_{\ell+1}}_{z_1,\ldots,z_{\ell+1}}
(\underline{p}_{\ell+1})\,=\, \sum_{p_{k,i}\in\CP^{(\ell+1)}}\,\,
\prod_{k=1}^{\ell+1} z_k^{\sum_i  p_{k,i}-\sum_i p_{k-1,i}}\,\,
\\ \times\frac{\prod\limits_{k=2}^{\ell}\prod\limits_{i=1}^{k-1}
(p_{k,i}-p_{k,i+1})_q!}
{\prod\limits_{k=1}^{\ell}\prod\limits_{i=1}^k
(p_{k+1,i}-p_{k,i})_q!\,\, (p_{k,i}-p_{k+1,i+1})_q!},\ee
and zero otherwise. Then,
 $\Psi^{\mathfrak{gl}_{\ell+1}}_{z_1,\ldots,z_{\ell+1}}(\underline{p}_{\ell+1})$
is a common solution of the eigenvalue problem (\ref{eiglat}).
\end{te}
Formula (\ref{main}) can be written also in the recursive form.
\begin{cor}
Let $\CP_{\ell+1,\ell}$
be a set of $\underline{p}_\ell=(p_{\ell,1},\ldots,p_{\ell,\ell})$
satisfying the conditions $p_{\ell+1,i}\geq p_{\ell,i}\geq
p_{\ell+1,i+1}$.  The following recursive relation holds:
\be\label{qtodarec}
 \Psi^{\mathfrak{gl}_{\ell+1}}_{z_1,\ldots,z_{\ell+1}}
(\underline{p}_{\ell+1})\,=\,
\sum_{\underline{p}_\ell\in\CP_{\ell+1,\ell}}\,\,
\Delta(\underline{p}_{\ell}) \,\,z_{\ell+1}^{\sum_i
p_{\ell+1,i}-\sum_i p_{\ell,i}}\,\,
Q_{\ell+1,\ell}(\underline{p}_{\ell+1},\underline{p}_{\ell}|q)
\Psi^{\mathfrak{gl}_{\ell}}_{z_1,\ldots,z_{\ell}}(\underline{p}_{\ell}),
\nonumber \ee where \be
Q_{\ell+1,\ell}(\underline{p}_{\ell+1},\underline{p}_{\ell}|q)\,=\,
\frac{1}{\prod\limits_{i=1}^{\ell} (p_{\ell+1,i}-p_{\ell,i})_q!\,\,
(p_{\ell,i}-p_{\ell+1,i+1})_q!},\\
\Delta(\underline{p}_{\ell})=\prod_{i=1}^{\ell-1}(p_{\ell,i}-p_{\ell,i+1})_q!\,\,\,. \ee
\end{cor}

\begin{rem}\label{gzgg} The representation \eqref{main} is a
 $q$-analog of  Givental's integral representation of the classical
$\mathfrak{gl}_{\ell+1}$-Whittaker function \cite{Gi2}, \cite{JK}:
\bqa\label{giv11}
\psi_{\underline{\lambda}}^{\mathfrak{gl}_{\ell+1}}
(x_1,\ldots,x_{\ell+1})\,=\,\int\limits_{\RR^{\frac{\ell(\ell+1)}{2}}}
\prod_{k=1}^{\ell}\prod_{i=1}^kdt_{k,i}\quad
e^{\frac{1}{\hbar}\mathcal{F}^{\mathfrak{gl}_{\ell+1}}(t) }, \eqa
where \bqa\label{intrep11}
\mathcal{F}^{{\mathfrak{gl}}_{\ell+1}}(t)=
\imath\sum\limits_{k=1}^{\ell+1} \lambda_k\Big(\sum\limits_{i=1}^{k}
t_{k,i}-\sum\limits_{i=1}^{k-1}t_{k-1,i}\Big)-
\sum\limits_{k=1}^{\ell} \sum\limits_{i=1}^{k}
\Big(e^{t_{k+1,i}-t_{k,i}}+e^{t_{k,i}-t_{k+1,i+1}}\Big), \nonumber
\eqa 
$\underline{\lambda}=(\lambda_1,\ldots,\lambda_{\ell+1})$,
$x_i:=t_{\ell+1,i},\,\,\,i=1,\ldots,\ell+1$ and $z_i=q^{\gamma_i}$, 
$\lambda_i=\gamma_i\log q$.  For the representation
theory derivation of this integral representation of
$\mathfrak{gl}_{\ell+1}$-Whittaker function see \cite{GKLO}.
The representation \eqref{main} of the $q$-Whittaker function  
turns into representation \eqref{giv11} the classical Whittaker
function in appropriate limit.

\end{rem}

As an example consider $\mathfrak{g}=\mathfrak{gl}_{2}$. Let 
$p_1:=p_{2,1}\in\ZZ$, $ p_2:=p_{2,2}\in\ZZ$ and $p:=p_{1,1}\in \ZZ$.
Then the  function \be\label{example} \Psi_{z_1,z_2}^{{\mathfrak gl}_2}
(p_{1},p_{2})\,\,=\,\,\sum_{p_{2}\leq p\leq p_{1}}\frac{ z_1^{p}
z_2^{p_{1}+p_{2}-p}} {(p_1-p)_q!(p-p_2)_q!},\hspace{1.5cm}
p_{1}\geq p_{2}\,, \ee
$$
\Psi_{z_1,z_2}^{{\mathfrak gl}_2}(p_{1},p_{2})=0, \hspace{3cm}
p_{1}<p_{2}\,,
$$
is a  solution of the system of equations:
\be
 \Big\{\,(1-q^{p_1-p_2+1})T_1\,+\,T_2\,\Big\}\,
 \Psi_{z_1,z_2}^{\mathfrak{gl}_2}(p_{1},p_{2})\,\,=\,\,(z_1+z_2)\,
 \Psi_{z_1,z_2}^{\mathfrak{gl}_2}(p_{1},p_{2}),\\
 T_1T_2\Psi_{z_1,z_2}^{{\mathfrak gl}_2}(p_{1},p_{2})\,\,=\,\,
 z_1z_2\,\Psi_{z_1,z_2}^{{\mathfrak gl}_2}(p_{1},p_{2}).\ee

Let us consider the following specialization of the $q$-deformed
$\mathfrak{gl}_{\ell+1}$-Whittaker function \be\label{spec}
\Psi^{\mathfrak{gl}_{\ell+1}}_{z_1,\ldots,z_{\ell+1}}(n,k)\,
:=\,\Psi^{\mathfrak{gl}_{\ell+1}}_{z_1,\ldots,z_{\ell+1}}( n+k,k,\ldots
,k). \ee

\begin{te}\label{Th12} 
$\Psi^{\mathfrak{gl}_{\ell+1}}_{z_1,\ldots,z_{\ell+1}}(n,k)$
satisfies  following  difference equation:
\be\label{eqdif}
\Big\{\prod_{i=1}^{\ell+1}(1-z_iT^{-1})\Big\}\,
\Psi^{\mathfrak{gl}_{\ell+1}}_{z_1,\ldots,z_{\ell+1}}(n,k)\,=\,
q^n \Psi^{\mathfrak{gl}_{\ell+1}}_{z_1,\ldots,z_{\ell+1}}(n,k),
 \ee where $T\cdot
f(n)=f(n+1)$.
\end{te}
\noindent {\it Proof}:  The proof is based on the explicit expression
(\ref{main}).  Let $\CP_{n,k}$ be a Gelfand-Zetlin pattern such that
 $(p_{\ell+1,1},\ldots,p_{\ell+1,\ell+1})=(n+k,k,\ldots,k)$.
Then, the relations $p_{\ell+1,i}\geq p_{\ell,i}\geq
p_{\ell+1,i+1}$ for  the elements of a Gelfand-Zetlin pattern imply 
 $p_{k,i\neq 1}=0$  and we have that \be
\Psi^{\mathfrak{gl}_{\ell+1}}_{z_1,\ldots,z_{\ell+1}}(n,k)
=\Big(\prod_{i=1}^{\ell+1}z_i^{k}\Big)
 \sum_{\CP_{n,k}} \frac{z_{\ell+1}^{n+k-p_{\ell,1}}}{(n+k-p_{\ell,1})_q!}
 \frac{z_{\ell}^{p_{\ell,1}-p_{\ell-1,1}}}{(p_{\ell,1}-p_{\ell-1,1})_q!}\cdots
 \frac{z_1^{p_{1,1}-k}}{(p_{{1,1}}-k)_q!}\\
 =\Big(\prod_{i=1}^{\ell+1}z_i^{k}\Big) \sum_{n_1+ \cdots +
 n_{\ell+1}=n}\frac{z_{\ell+1}^{n_{\ell+1}}}{(n_{\ell+1})_q!}
 \cdots\frac{z_1^{n_1}}{(n_1)_q!}.
\ee 
Introduce the generating function
$$
\Psi^{\mathfrak{gl}_{\ell+1}}_{z_1,\ldots,z_{\ell+1}}(t,k)
=\sum_{n\in\IZ}\,t^n\,\Psi^{\mathfrak{gl}_{\ell+1}}_{z_1,\ldots,z_{\ell+1}}(n,k)
=\prod_{j=1}^{\ell+1}\frac{z_i^k}{\prod_{m=0}^{\infty}(1-tz_iq^m)},
$$
where we use the identity
$$
\frac{1}{\prod_{m=0}^{\infty}(1-xq^m)}=\sum_{n=0}^{\infty}\,\frac{x^n}{(n)_q!}.
$$
Due to the fact that
$\Psi^{\mathfrak{gl}_{\ell+1}}_{z_1,\ldots,z_{\ell+1}}(n,k)=0$
 for $n<0$, the generating function
 $\Psi^{\mathfrak{gl}_{\ell+1}}_{z_1,\ldots,z_{\ell+1}}(t,k)$ is
 regular at $t=0$.
It is easy to check now the following identity
$$
\prod_{j=1}^{\ell+1}(1-tz_i)\,\Psi^{\mathfrak{gl}_{\ell+1}}_{z_1,\ldots,z_{\ell+1}}(t,k)
=\Psi^{\mathfrak{gl}_{\ell+1}}_{z_1,\ldots,z_{\ell+1}}(qt,k).
$$
Expanding the latter relation in powers of $t$, we obtain
\eqref{eqdif} for the coefficients of 
$\Psi^{\mathfrak{gl}_{\ell+1}}_{z_1,\ldots,z_{\ell+1}}(t,k)$  $\Box$

\begin{rem} The difference equation \eqref{eqdif}
for the specialized $q$-Whittaker function
$\Psi^{\mathfrak{gl}_{\ell+1}}_{z_1,\ldots,z_{\ell+1}}(n,k)$
can be derived   directly from the system of equations \eqref{eiglat}
for the  non-specialized   $q$-deformed Whittaker function
$\Psi^{\mathfrak{gl}_{\ell+1}}_{z_1,\ldots,z_{\ell+1}}(p_1,p_2,\ldots,p_{\ell+1})$
and the condition
$$
\Psi^{\mathfrak{gl}_{\ell+1}}_{z_1,\ldots,z_{\ell+1}}(p_1,p_2,\ldots,p_{\ell+1})=0
$$
outside dominant domain $p_1\geq\ldots \geq p_{\ell+1}$.
\end{rem}

 \begin{lem}\label{lemma}
 The following integral representation for the specialized  $q$-deformed
$\mathfrak{gl}_{\ell+1}$-Whittaker functions holds
 \be\label{qWhit}
  \Psi^{{\mathfrak gl}_{\ell+1}}_{\underline{z}}(n,k)\,=\,
\Big(\prod_{i=1}^{\ell+1}z_i^{k}\Big)\, 

\oint\limits_{t=0}\,\!\!\frac{dt}{2\pi\imath\,t}\,\,
   t^{-n}\prod_{i=1}^{\ell+1}\,\Gamma_q(z_it),
 \ee
where
$$
 \Gamma_q(x)\,=\,\prod_{n=0}^{\infty}\frac{1}{1-q^nx}.
$$
\end{lem}

\noindent {\it Proof:} Using the identity
$$
 \prod_{n=0}^\infty\frac{1}{1-xq^n}\,=\,
 \sum_{m=0}^{\infty}\frac{x^m}{(m)!_q},
$$
one obtains, for $n\geq 0$, that 
 \bqa \Psi^{{\mathfrak gl}_{\ell+1}}_{\underline{z}}(n,k)\,=\,
\Big(\prod_{i=1}^{\ell+1}z_i^{k}\Big)  \oint\limits_{t=0}\!\!\frac{dt}{2\pi\imath\,t}
  \,t^{-n}\,\,\prod_{i=1}^{\ell+1}\prod_{m=0}^\infty
  \frac{1}{1-z_itq^m}\\ \nonumber
  =\,\Big(\prod_{i=1}^{\ell+1}z_i^{k}\Big) \sum_{n_1+\ldots+n_{\ell+1}=n}
  \frac{z_1^{n_1}}{(n_1)_q!}\cdot\ldots\cdot
  \frac{z_{\ell+1}^{n_{\ell+1}}}{(n_{\ell+1})_q!}\,
 \eqa
For $n<0$, we obviously have that 
$\Psi^{{\mathfrak gl}_{\ell+1}}_{\underline{z}}(n,k)=0\,\,\,$
$\Box$

The corresponding integral representation 
for the classical $\mathfrak{gl}_2$-Whittaker function is  given
 by  the Mellin-Barnes
representation for  the $\mathfrak{gl}_2$-Whittaker function
\be\label{mellin2}
\psi^{\mathfrak{gl}_2}_{\lambda_1,\lambda_2}(x_1,x_2)\,=\,
e^{\frac{\i}{\hbar}(\lambda_1+\lambda_2)x_2}\,\int_{i\sigma-\infty}^{\i\sigma+\infty}
\!\!d\lambda\,\, e^{\frac{\i}{\hbar}\lambda(x_1-x_2)}\,\,
\Gamma\Big(\frac{\lambda-\lambda_1}{\imath\hbar}\Big)\,
\Gamma\Big(\frac{\lambda-\lambda_2}{\imath\hbar}\Big), \ee
 where $\sigma > \max\{Im \lambda_j,\,j=1,\ldots,\ell+1\}$.
 
 \begin{rem} The expression
\be\label{redmain}
\Psi^{\mathfrak{gl}_{\ell+1}}_{z_1,\ldots,z_{\ell+1}}(n,k)\,=\,
 \Big(\prod_{i=1}^{\ell+1}z_i^{k}\Big) \sum_{n_1+\ldots+n_{\ell+1}=n}\,\,
 \frac{z_1^{n_1}}{(n_1)_q!}\cdot\ldots\cdot
 \frac{z_{\ell+1}^{n_{\ell+1}}}{(n_{\ell+1})_q!}\,, \qquad n\geq 0,\\
\Psi^{\mathfrak{gl}_{\ell+1}}_{z_1,\ldots,z_{\ell+1}}(n,k)\,=\,0,\qquad \,n<0
\ee 
is a $q$-analog of the 
Givental integral representation for the equivariant
Gromov-Witten invariants of $X=\IP^{\ell}$ {\rm \cite{Gi3}}
\be\label{redgiv11}
f_{\underline{\lambda}}(T)\,=\,\int\limits_{\RR^{\ell}}\prod_{k=1}^\ell
dt_{k,1}\,\ \, e^{\frac{1}{\hbar}\mathcal{F}(t)},
 \ee
 where $\underline{\lambda}=(\lambda_1,\ldots,\lambda_{\ell+1})$,
 $T:=t_{\ell+1,1}$ and \bqa\nonumber
\mathcal{F}(t)\,=\,\imath\lambda_1t_{11}+
\sum\limits_{k=1}^\ell\imath\lambda_{k+1} (t_{k+1,1}-t_{k,1})\,-\,
e^{t_{11}}-\sum\limits_{k=1}^\ell e^{t_{k+1,1}-t_{k,1}}.\eqa
The representation \eqref{redmain} for specialized $q$-Whittaker
function  turns into \eqref{redgiv11} in appropriate limit. 
\end{rem}

\section{Counting  holomorphic sections}  

In this Section we are going to provide an interpretation of the  explicit
expressions for $q$-deformed  
specialized class one  $\mathfrak{gl}_{\ell+1}$-Whittaker functions in terms
of traces of operators acting on the spaces 
of holomorphic sections of line bundles 
on infinite-dimensional  manifolds. For this aim, we first 
consider an auxiliary problem of counting  holomorphic sections on 
finite-dimensional manifolds approximating  the infinite-dimensional
ones. The relevant finite-dimensional manifolds are spaces of
the quasi-maps of $\IP^1$ to $GL_{\ell+1}(\IC)$-homogeneous 
spaces. 

\subsection{Space of quasi-maps}

Let us start with recalling the  general construction of the  quasi-map
compactification of the space of holomorphic maps of $\IP^1$  to
the partial flag spaces of complex Lie group $GL_{\ell+1}$ 
due to Drinfeld. Let $\alpha_i$,
$i=1,\ldots ,\ell$,  be a set of simple roots of the complex Lie algebra
$\mathfrak{gl}_{\ell+1}$. To any  ordered subset of simple roots
$\{\alpha_{i_1},\ldots,\alpha_{i_r}\}$ indexed by an ordered subset
$I^P=\{i_1<\ldots<i_r\}\subset\{1,\ldots,\ell\}$ one can associate a
parabolic subgroup $P\subset GL_{\ell+1}$.  Namely, let
$B\subset GL_{\ell+1}$ be the subgroup of upper-triangular
matrices generated by Cartan torus and 
one-parameter unipotent subgroups
corresponding to positive simple roots. Then a parabolic subgroup
$P$ is generated by $B$ and one-parameter unipotent subgroups
corresponding to negative roots $-\alpha_i$ such that $i\notin I^P$.
In particular, when $r=\ell$ one gets $P=B$, and the corresponding
homogeneous space $GL_{\ell+1}/B$ coincides with the full flag
space.  On the other hand for a parabolic subgroup $P_0\subset
GL_{\ell+1}$ associated to the first simple root (i.e.
$I^{P_0}=\{1\}\subset\{1,2,\ldots,\ell\}$), the corresponding
homogeneous space $GL_{\ell+1}/P_0$ is isomorphic to the
projective space $\IP^\ell$.  Partial flag 
spaces $GL_{\ell+1}/P$ possess canonical  projective embeddings
\be\label{Pluck} \pi:\,GL_{\ell+1}/P\to \Pi=\prod_{j\in I^P}\,
\IP^{n_j-1},\hspace{2cm} n_j=(\ell+1)!/j!\,\,(\ell+1-j)!\,. \ee

The group $H^2(GL_{\ell+1}/P,\IZ)=\IZ^r$ is naturally isomorphic to a
sublattice of the weight lattice of $\mathfrak{sl}_{\ell+1}$ and is
spanned by the weights $\omega_i$ indexed by $I^P$.
 Let $\CL_j$, $j=1,\ldots, r$, be the line bundles  on
$GL_{\ell+1}/P$ obtained  
as pull backs of $\CO(1)$ form the direct factors $\IP^{n_j-1}$ in the
right hand side (r.h.s.)  of \eqref{Pluck}. 
The lattice $H^2(GL_{\ell+1}/P,\IZ)=\IZ^r$ is
generated by the first Chern classes $c_1(\CL_i)$.

Let $\CM_{\underline{d}}(GL_{\ell+1}/P)$  be a non-compact space of
 holomorphic  maps of $\IP^1$ of multi-degree $\underline{d}\in H^2(GL_{\ell+1}/P,\IZ)$
 to the flag space $GL_{\ell+1}/P$. Due to \eqref{Pluck}, 
$\CM_{\underline{d}}(GL_{\ell+1}/P)$
is a subspace of the product of space
$\CM_{d_j}(\IP^{n_j-1})$.  Explicitly, 
 each $\CM_{d_j}(\IP^{n_j-1})$ can be described as a set of
collections of $n$, {\it relatively prime}
 polynomials of degree $d$, up to a common constant
factor. The space $\CM_{d_j}(\IP^{n_j-1})$ allows for a compactification
by the  space of quasi-maps $\CQ\CM_{d_j}(\IP^{n_j-1})$ defined as a set of
collections of $n$, polynomials of degree $d$, up to a common
constant factor. The space  of quasi-maps
$\CQ\CM_{\underline{d}}(GL_{\ell+1}/P)$
is then constructed  as a closure of  $\CM_{\underline{d}}(GL_{\ell+1}/P)$
in $\prod_j \CQ\CM_{\underline{d}}(\IP^{n_j})$. Thus defined
$\CQ\CM_{\underline{d}}(GL_{\ell+1}/P)$ is (in general singular)
irreducible projective variety. A 
small resolution of this space is known due to  \cite{La}, \cite{Ku}.

On the space of holomorphic maps $\CM_{\underline{d}}(GL_{\ell+1}/P)$
of $\IP^1$ to  $GL_{\ell+1}/P$, 
there is a natural action of the group $\IC^*\times GL_{\ell+1}$ (and, thus,
of its maximal compact subgroup $S^1\times U_{\ell+1}$). Here, 
the  action of
$GL_{\ell+1}$ is induced by the standard action on flag spaces and the
action of $\IC^*$ is induced by the action  of 
$\IC^*$ on $\IP^1$  given by  $(y_1,y_2)\to (\xi
y_1,y_2)$ in homogeneous coordinates $(y_1,y_2)$ on $\IP^1$.
This action of $\IC^*\times GL_{\ell+1}$ can be
extended to an action on the space $\CQ\CM_{\underline{d}}(GL_{\ell+1}/P)$ of
quasi-maps.

In the following we consider a parabolic subgroup $P_0\subset GL_{\ell+1}$ associated
to the first simple root (and thus
$I^{P_0}=\{1\}\subset\{1,2,\ldots,\ell\}$. The corresponding
homogeneous space $GL_{\ell+1}/P_0$ is a projective space
$\IP^\ell$. The space of quasi-maps
\footnote{The compactification of $\CM_d(\IP^{\ell})$ 
by the space $\CQ\CM_d(\IP^{\ell})$ of quasi-maps 
arises naturally in the linear sigma-model
  description of Gromov-Witten invariants of projective spaces \cite{W}.}  
 $\QM_d(\IP^\ell)$ is a non-singular
projective variety $\IP^{(\ell+1)(d+1)-1}$. A quasi-map
$\phi \in \CQ\CM_d(\IP^{\ell})$  is given by a collection
$$
\bigl(a_0({y}):a_1({y}):\,\ldots\, a_\ell({y})\bigr),
$$
of homogeneous polynomials $a_i({y})$ in variables
${y}=(y_1,y_2)$ of degree $d$ 
$$
a_k({y})\,=\,\sum_{i=0}^da_{k,j}\,y_1^{j}y_2^{d-j},
\hspace{2cm}k=0,\ldots,\ell.
$$
considered up to the multiplication of all $a_i(y)$'s by a nonzero
complex number.  The action of $(\xi,g)\in \IC^*\times
GL_{\ell+1}$ on $\QM_d(\IP^\ell)$ is given by
 \be\xi\,:\quad\bigl(a_0({y}):a_1({y}):\,\ldots\,
a_\ell({y})\bigr)\,\longmapsto\,\bigl(
a_0({y}^\xi):a_1({y}^\xi):\,\ldots:\,
a_\ell({y}^\xi)\bigr),\\
g\,:\quad\bigl(a_0({y}):a_1({y}):\,\ldots\,:
a_\ell({y})\bigr)\longmapsto\Big(\,
\sum_{k=1}^{\ell+1}g_{1,k}a_{k-1}({y}):\,\ldots\,:
\sum_{k=1}^{\ell+1}g_{\ell+1,k}a_{k-1}({y})\,\Big),\ee where
 $g=\|g_{ij}\|$ and ${y}^\xi=(\xi y_1,y_2)$.

\subsection{Generating functions of holomorphic sections}

Let $\CO(1)$ be a standard line bundle on $\IP^{(\ell+1)(d+1)-1}$. 
The space of sections of the line bundle $\CO(n):=\CO(1)^{\otimes n}$ on
$\CQ\CM_d(\IP^{\ell})$ is naturally a $\IC^*\times
GL_{\ell+1}$-module. We are interested in calculating the
corresponding character. 

Let $T\in GL_{\ell+1}$ be a Cartan torus,
$H_1,\ldots,H_{\ell+1}$ be a basis in $\Lie(T)$, and $L_0$ be a
generator of $\Lie(\IC^*)$. The equivariant cohomology of a point with
respect to the maximal compact subgroup $G=S^1\times U_{\ell+1}$ of
$\IC^*\times GL_{\ell+1}$ can be described as
$$
 H^*_G({\rm pt},\C)\,=\,
 \C[\la_1,\ldots,\la_{\ell+1}]^{\mathfrak{S_{\ell+1}}}\otimes\C[\hbar],
$$
where $\la_1,\ldots,\la_{\ell+1}$ and $\hbar$ are associated with the
generators $H_1,\ldots,H_{\ell+1}$ and $L_0$ respectively.

Let $\CL_k$ be a one-dimensional $GL_{\ell+1}$-module such that
$H_i\CL_k=k\CL_k$, for $i=1,\ldots,\ell+1$. 
Cohomology groups $H^*(\QM_d(\IP^\ell),\CO(n))\otimes \CL_k$
have a natural structure of $\IC^*\times GL_{\ell+1}(\IC)$-module.
We denote by $\CL_k(n)=\CO(n)\otimes \CL_k$  the 
line bundles $\CO(n)$ twisted by one-dimensional $GL_{\ell+1}$-module
$\CL_k$. 

Let $A^{(d)}_{n,k}(\underline{z},q)$,  be the
character of the  $\IC^*\times GL_{\ell+1}$-module
$\CV_{n,k,d}\,=\,H^0(\QM_d(\IP^\ell),\CL_k(n))$ , 
$n\geq 0$,
$$
 A^{(d)}_{n,k}(\underline{z},q)\,=\,
 \Tr_{\CV_{n,k,d}}\,q^{L_0}\,e^{\sum\lambda_iH_i},
$$
where we assume that $q\in \IC^*, |q|<1$. 
This character can be straightforwardly calculated as follows.
The space $\CV_{n,k,d}=H^0(\QM_d(\IP^\ell),\CL_k(n))$ can be identified
with the space of degree $n$  homogeneous polynomials in
$(\ell+1)(d+1)$ variables $a_{k,i}$, for $k=0,\ldots,\ell$ and
$i=0,\ldots,d$. Define 
$$\CV_{k,d}=\oplus_{n=0}^{\infty} \CV_{n,k,d},$$
and the grading on $\CV_{k,d}$ is defined 
by the eigenvalue decomposition with respect to the action of an operator $D$
$$
t^D:\,\CV_{n,k,d}\rightarrow t^n\,\,\CV_{n,k,d},\qquad t\in \IC^*.
$$
The action of the subgroup $(\IC^*\times T)\subset G(\IC)=\IC^*\times
GL_{\ell+1}$ is given by \be e^{\sum\lambda_iH_i}\,:\quad
\bigl(a_0({y}):a_1({y}):\,\ldots\,
a_\ell({y})\bigr)\,\longmapsto\,\bigl(
e^{\lambda_1}a_0({y}):e^{\lambda_2}a_1({y}):\,
\ldots:\,e^{\lambda_{\ell+1}}a_\ell({y})\bigr),\ee 
where 
$$
a_k({y})\,=\,\sum_{j=0}^da_{k,j}\,y_1^{j}y_2^{d-j},
\hspace{2cm}k=0,\ldots,\ell.
$$
The action of the generator $L_0$ of $\IC^*$ is as follows \be
q^{L_0}\,:\quad a_{k,j}\,\longmapsto\, q^j\,\,a_{k,j}. \ee

\begin{prop}\label{refmainth}  For the  $\IC^*\times GL_{\ell+1}$-character of 
the module $\CV_{n,k,d}$, the following integral representation holds
\be\label{finited}   
A^{(d)}_{n,k}(\underline{z},q)\,=\,
 \Tr_{_{\CV_{n,k,d}}}\,q^{L_0}\,e^{\sum\lambda_iH_i}\,=\,
 \Big(\prod_{i=1}^{\ell+1}z_i^k\Big)\,\oint_{t=0}\,
 \frac{dt}{2\pi\imath\,t^{n+1}}\,
 \prod_{m=1}^{\ell+1}\prod_{j=0}^d\frac{1}{(1-tq^jz_m)},
\ee where $\underline{z}=(z_1,\ldots,z_{\ell+1})$ and
$z_m=e^{\lambda_m}$.
\end{prop}

\noindent {\it Proof}: A simple calculation  gives us that 
\be
 A^{(d)}_k(\underline{z},t,q)\,=\,
 \Tr_{_{\CV_{k,d}}}\,t^{D}\,q^{L_0}\,e^{\sum\lambda_iH_i}\,=\,
 \Big(\prod_{i=1}^{\ell+1}z_i^k\Big)
 \prod_{m=1}^{\ell+1}\prod_{j=0}^d\frac{1}{(1-tq^jz_m)}.
\ee
The projection on the subspace of $\CV_{k,d}$ of the grading $n$ with
respect to $D$ can be realized by taking a residue,  
\be
 \Tr_{_{\CV_{n,k,d}}}\,q^{L_0}\,e^{\sum\lambda_iH_i}\,=\,
\oint_{t=0}\,  \frac{dt}{2\pi\imath\,t^{n+1}}\,\,\,
 \Tr_{_{\CV_{k,d}}}\,t^{D}\,q^{L_0}\,e^{\sum\lambda_iH_i}.
\ee
This gives us the integral expression \eqref{finited} $\Box$

\subsection{Equivariant Euler characteristic of line bundles on 
$\CQ\CM_d(\IP^{\ell})$}\label{modelcoh}

Characters \eqref{finited} of the space of holomorphic sections  can be
related to equivariant holomorphic Euler characteristics  of line bundles
on $\CQ\CM_d(\IP^{\ell})$. First we recall the standard facts about
line bundles on projective spaces. Line bundles $\CO(n)$ on projective spaces $\IP^N$ 
are equivariant with respect to the 
standard action of $U_{N+1}$ on $\IP^N$. 
The $U_{N+1}$-equivariant Euler characteristic of $\CO(n)$ is given by 
the  character 
\be\label{RRHproj} 
\chi_{_{U_{N+1}}}(\IP^N,\CO(n))\,=\,\sum_{m=0}^{N}\,
(-1)^m\,\Tr_{H^m(\IP^N,\CO(n))}\,\,
 e^{\sum\lambda_iH_i}\,,\qquad  e^{\sum\lambda_iH_i}\in U_{N+1}.\ee
Cohomology groups of  $\CO(n)$ on projective space $\IP^N$ 
have the following properties (see e.g. \cite{OSS}) \bqa\label{oandl} \dim
H^m(\IP^N,\CO(n))&=&0, \qquad m\neq 0,N,\nonumber\\ 
\dim H^N(\IP^N,\CO(n))&=&0,\qquad n\geq 0,\\
 \dim H^{0}(\IP^N,\CO(n))&=&0,\qquad n<0.\nonumber
\eqa  
Taking into account \eqref{oandl} the expression  \eqref{RRHproj} reduces to 
\be\label{decomptwo}
\chi_{U_{N+1}}(\IP^N,\CO(n))\,=\,\Tr_{H^0(\IP^N,\CO(n))}\,\,
  e^{\sum\lambda_iH_i}\,+(-1)^N\,\Tr_{H^N(\IP^N,\CO(n))}\,\,
  e^{\sum\lambda_iH_i}\,.\ee
We have $\QM_d(\IP^\ell)=\IP^{(\ell+1)(d+1)-1}$ and, thus, for $n\geq
0$,  we can identify 
$A^{(d)}_{n,k}(\underline{z},q)$ with
equivariant holomorphic Euler characteristic of 
$\CL_k(n)$  
\be \nonumber
A^{(d)}_{n,k}(\underline{z},q)=
\Tr_{H^0(\QM_d(\IP^\ell),\CL_k(n))}\,\,
 e^{\hbar L_0+\sum\lambda_iH_i}=
\chi_G(\QM_d(\IP^\ell),\CL_k(n)),\qquad n\geq 0,\ee
where $G=S^1\times U_{\ell+1}$.
The equivariant Euler characteristic of a holomorphic vector
bundle on the projective space possesses  a canonical holomorphic integral
representation. According to the Riemann-Roch-Hirzebruch (RRH) theorem,
one can express the $U_{N+1}$-equivariant holomorphic Euler
characteristic of a $U_{N+1}$-equivariant vector bundle $\CE$ on
$\IP^N$ as follows 
\be\label{RRHzeroone}
 \chi_{U_{N+1}}(\IP^N,\CE)=
\sum_{m=0}^{N}\,(-1)^m \Tr_{H^m(\IP^N,\CE)}\,\,
 e^{\sum\lambda_iH_i}=\\
 =\,\<{\rm Ch}_{U_{N+1}}(\CE)\,\,{\rm Td}_{U_{N+1}}
 (\CT\IP^N)\,,[\IP^N]\>,
\ee  where $H_i$ are generators of the  Lie algebra
$\mathfrak{gl}_{N+1}$, $\CT\IP^N$ is the tangent bundle to $\IP^N$, 
  ${\rm Ch}_{U_{N+1}}(\CE)$ is a
$U_{N+1}$-equivariant Chern character of $\CE$ and  ${\rm
Td}_{U_{N+1}}(\CE)$ is a $U_{N+1}$-equivariant Todd genus of $\CE$
\cite{A}.

The tangent bundle $\CT\IP^N$ to the 
projective space $\IP^N$ is $U_{N+1}$-equivariantly
 stable-equivalent to $\CO(1)^{\oplus
  (N+1)}$ as the following lemma shows. 

\begin{lem}\label{EEseq}
The following relation holds in $U_{N+1}$-equivariant  topological
$K$-theory $K_{U_{N+1}}(\IP^{N})$
\be\label{Eseq}
[\CT\IP^{N}]\oplus [\CO]=[\CO(1)]^{\oplus (N+1)},
\ee
\end{lem}
where $[\CE]$ is a class of a vector bundle $\CE$ in $K_{U_{N+1}}(\IP^N)$.

\noindent {\it Proof}:  For a tangent sheaf to the complex projective
space $\IP^N$ we have the Euler exact sequence (see e.g. \cite{GH})
\be\label{exactseqproj}
 0\longrightarrow\CO\longrightarrow\CO(1)^{\oplus(N+1)}
 \longrightarrow \CT\IP^N\longrightarrow 0
\ee
The maps \eqref{exactseqproj} are explicitly  $U_{N+1}$-equivariant
and, thus, we obtain the relation \eqref{Eseq} in
$U_{N+1}$-equivariant $K$-groups of
$\IP^N$ $\Box$

Lemma \ref{EEseq} and the fact that the Todd class
depends only on stable equivalence class of a vector bundle
allows us to rewrite RRH-theory on projective spaces
as follows
\be\label{RRHzero}
\chi_{_{U_{N+1}}}(\IP^N,\CE)=
\<{\rm Ch}_{U_{N+1}}(\CE)\,\,{\rm Td}_{U_{N+1}}
 (\CO(1)^{\oplus(N+1)}\,,[\IP^N]\>.
\ee 
In the following we will consider only the case of line bundles
and thus we take $\CE=\CO(n)$, $n\in \IZ$. The
pairing of the cohomology classes with the fundamental class entering
the formulation of RRH-theorem can be expressed explicitly using a
particular model for the cohomology ring $H^*(\IP^N,\IC)$. The
cohomology ring $H^*(\IP^N,\IC)$ is generated by an element $x\in
H^2(\IP^N,\IC)$ with a single relation $x^{N+1}=0$ \be\label{cohpt}
H^*(\IP^{N},\IC)=\IC[x]/x^{N+1}. \ee The $U_{N+1}$-equivariant analog
of \eqref{cohpt} is given by
$$
H_{U_{N+1}}^*(\IP^{\ell},\IC)=\IC[x]\otimes \IC[\lambda_1,\cdots ,
\lambda_{N+1}]^{\mathfrak{S}_{N+1}}\Big/\Big(
\prod_{j=1}^{N+1}(x-\lambda_j)\Big),
$$
which  is naturally a module over
$H_{U_{N+1}}^*(\pt,\IC)=\IC[\lambda_1,\cdots,
\lambda_{N+1}]^{\mathfrak{S}_{N+1}}$ where
$\mathfrak{S}_{N+1}$ is the permutation group of a set of $N+1$
elements. The pairing of an element of
$H^*_{U_{N+1}}(\IP^N,\IC)$ represented by $P(x,\lambda)$
with a $U_{N+1}$-equivariant fundamental
cycle $[\IP^N]$ can be expressed in terms of the integral
 $$
\<P(\lambda),[\IP^N]\>=\frac{1}{2\pi\imath}\, \oint_{C_0}
dx\,\frac{P(x,\lambda)}{\prod_{ j=1}^{N+1}(x-\lambda_j)},
$$
where the integration contour $C_0$ encircles the poles
$x=\lambda_j$, $j=1,\ldots, (N+1)$.
 The pairing for  $H^*(\IP^N,\IC)$ is obtained by a
specialization $\lambda_j=0$, $j=1,\cdots, (N+1)$. The equivariant
Chern character and Todd class can be written in terms of this model
of $H^*_{U_{N+1}}(\IP^N,\IC)$ as follows (see e.g. \cite{H})
$$
 {\rm Ch}_{U_{N+1}}(\CO(n))\,=\,e^{nx},\hspace{1.5cm}
 {\rm Td}_{U_{N+1}}\bigl(\CO(1)^{\oplus(N+1)}\bigr)\,=\,
 \prod_{j=1}^{N+1}\frac{(x-\lambda_j)}{1-e^{-(x-\lambda_j)}}.
$$
Therefore we have the following integral representation of the
equivariant holomorphic Euler characteristic
 ($t=e^{-x}\,,\,z_i=e^{\lambda_i},\,i=1,\ldots,N+1$):
\be\label{RRHqzero}
 \chi_{U_{N+1}}(\underline{z})\,=\, \<{\rm
 Ch}_{U_{N+1}}(\CO(n))\,{\rm
 Td}_{U_{N+1}}\bigl(\CO(1)^{\oplus(N+1)}\bigr),[\IP^N]\>\\
 =\,\frac{1}{2\pi \imath}\,\oint_{C_0} \,
 \frac{dx}{\prod_{i=1}^{N+1}(x-\lambda_i)}\,\,
 e^{nx}\,\,
 \prod_{i=1}^{N+1}\frac{(x-\lambda_i)}{(1-e^{-(x-\lambda_i)})}\ee
\be\label{Hcharone}=\,-\frac{1}{2\pi \imath}\,\oint_{C_0}
\,\frac{dt}{t^{n+1}}\,\prod_{i=1}^{N+1}\frac{1}{1-tz_i}\,,\ee where
in the last expression the integration contour $C_0$ encircles the
poles $t=z_j^{-1}, j=1,\ldots,\ell+1$. The integral representation
\eqref{Hcharone} can be obtained directly using a particular
realization of ($U_{N+1}$-equivariant) $K$-theory on $\IP^N$ (see
e.g. \cite{A}). The $K$-group $K(\IP^N)$ is generated by a
class $t$ of the line bundle $\CO(1)$ satisfying  the relation
$(1-t)^{N+1}=0$. We have the following isomorphisms for
($U_{N+1}$-equivariant) $K$-groups of projective spaces
\be\label{Ktheory}
K(\IP^N)=\IC[t,t^{-1}]/(1-t)^{N+1},\qquad
K_{U_{N+1}}(\IP^N)= \IC[t,t^{-1},z,z^{-1}]\,\Big
/\prod_{j=1}^{N+1}(1-tz_j)\,.
\ee
The equivariant analog of the pairing with the fundamental class of
$\IP^N$ in $K$-theory  is given by
\be\label{pairingK}
\<R\,,[\IP^N]\>_K\,=\,-\frac{1}{2\pi \imath}\,\oint\limits_{C_0}
\, \frac{dt}{t} \,\,\frac{R(t)}{\prod_{j=1}^{N+1}(1-tz_j)},
\ee
where $R(t)$ is a rational function representing an element
of $K_{U_{N+1}}(\IP^N)$ and the  integration contour $C_0$ encircles
the poles $t=z_j^{-1}$, $j=1,\ldots ,(N+1)$.

Using the representation of the pairing \eqref{pairingK}
one can represent RRH expression for the Euler characteristic as
\be
 \chi_{U_{N+1}}(\underline{z})\,=\,\<[\CO(n)],[\IP^N]\>_K=
-\frac{1}{2\pi \imath}\,\oint_{C_0}
\,\frac{dt}{t^{n+1}}\,\prod_{i=1}^{N+1}\frac{1}{1-tz_i}\,.\ee
This  reproduces the representation \eqref{Hcharone}.

Now we would like to apply the integral representation for equivariant
Euler  characteristic to the $S^1\times U_{\ell+1}$-equivariant line bundle
$\CL_k(n)$ on $\CQ\CM_d(\IP^{\ell})$.

Consider the $S^1\times U_{\ell+1}$-equivariant cohomology of the
projective space $\IP(V_{(\ell+1)(d+1)})$ where the vector space
$V_{(\ell+1)(d+1)}=\oplus_{j=1}^{\ell+1}\,\oplus_{m=0}^{d}
V_{j,m}$
 has the  structure of an $S^1$-module with an  action given by
$$
e^{\imath\theta}:\,V_{j,m}\rightarrow e^{\imath m \theta}\,V_{j,m},\qquad
\dim V_{j,m}=1,\quad \theta \in S^1,
$$
and each $V_{m}=V_{1,m}\oplus V_{2,m}\oplus \ldots V_{\ell+1,m}$
is standard $U_{\ell+1}$-module.  Then for the $G=S^1\times U_{\ell+1}$-equivariant
cohomology of $\IP(V_{(\ell+1)(d+1)})$ we have an isomorphism
\be\label{modelfor}
H_{S^1\times U_{\ell+1}}^*(\IP(V_{(\ell+1)(d+1)}))\,=\,
\IC[x,\lambda,\hbar]\Big/\prod_{j=1}^{\ell+1}\prod_{m=0}^{d}(x-\lambda_j-\hbar m),
\ee
where $x$ is a generator of $H^*(\IP(V_{(\ell+1)(d+1)}),\IC)$. The
pairing with the $S^1\times U_{\ell+1}$-equivariant 
fundamental cycle $[\IP(V_{(\ell+1)(d+1)})]$ can be expressed in
the form of the contour integral \be\label{pairing}
 \<P(\lambda,\hbar),[\IP(V_{(\ell+1)(d+1)})]\>\,=\,
\frac{1}{2\pi\imath}\oint_C
 \frac{P(x,\lambda,\hbar)\,\,dx}{\prod_{ j=1}^{\ell+1}
 \prod_{m=0}^{d}
 (x-\lambda_j - \hbar m)},
\ee where the integration contour
$C$ encircles the poles $x=\lambda_j+m\hbar$, $j=1,\ldots ,\ell+1$,\
 $m=0,1,\ldots ,d$.

Specializing to the  action of $G=S^1\times U_{\ell+1}$ on
$\QM_d(\IP^{\ell})\simeq\IP^{(\ell+1)(d+1)-1}$ described in Section 2.1 we
obtain
$$
 {\rm Ch}_{G}(\CL_k(n))\,=\,
 e^{nx+k(\lambda_1+\ldots+\lambda_{\ell+1})},\hspace{1.5cm}
 {\rm Td}_G (\CT\IP^{(\ell+1)(d+1)-1})\,=\,
 \prod_{i=1}^{\ell+1}\prod_{m=0}^d
 \frac{x- m\hbar-\lambda_i}
 {1-e^{\lambda_i+ m\hbar-x}}.
$$
Let $q=e^\hbar$, $t=e^{-x}$, and $z_i=e^{\lambda_i},\,1\leq
i\leq\ell+1$, then
\be\label{MainChar}
 \chi_G(\QM_d(\IP^\ell),\CL_k(n))\,=\,
 \bigl\<{\rm Ch}_{G}(\CL_k(n))\,
 {\rm Td}_G (\CT\IP^{(\ell+1)(d+1)-1}),
 [\IP^{\ell}]\bigr\>\\
 =\,\frac{1}{2\pi \imath}\,\oint_C dx\,\,
 \prod_{i=1}^{\ell+1}\prod_{m=0}^{d}\frac{1}{(x-\lambda_i-m\hbar)}\,\,
 e^{nx+k(\lambda_1+\ldots+\lambda_{\ell+1})}\,\,
 \prod_{i=1}^{\ell+1}\prod_{n=0}^{d}
 \frac{(x-\lambda_i-m\hbar)}{(1-e^{-(x-\lambda_i-m\hbar)})}\\
 =\,-\Big(\prod_{i=1}^{\ell+1}z_i^k\Big)\,\oint\limits_{C}\,
 \frac{dt}{2\pi\imath\,t^{n+1}}\,
 \prod_{i=1}^{\ell+1}\prod_{m=0}^d\frac{1}{(1-tz_iq^m)}.
\ee
For $n\geq 0$ one has the identity
$$
\chi_G(\QM_d(\IP^\ell),\CL_k(n))
=\Tr_{\CV_{n,k,d}}\,q^{L_0}\,e^{\sum\lambda_iH_i}.
$$
Deforming  the contour for $n\geq 0$
 we obtain  the following integral representation for the character
$$
\Tr_{\CV_{n,k,d}}\,q^{L_0}\,e^{\sum\lambda_iH_i}\,=\,
\,\Big(\prod_{i=1}^{\ell+1}z_i^k\Big)\,\oint\limits_{t=0}\,
 \frac{dt}{2\pi\imath\,t^{n+1}}\,
 \prod_{i=1}^{\ell+1}\prod_{m=0}^d\frac{1}{(1-tz_iq^m)}
$$
which coincides  with \eqref{finited}.

\begin{rem} Without a restriction $n\geq 0$, the integral representation
for the equivariant Euler characteristic can be represented  as a
difference of two terms
\be\nonumber
 \chi_G(\QM_d(\IP^\ell),\CL_k(n))\,=\,\\
=\Big(\prod_{i=1}^{\ell+1}z_i^k\Big)\,\oint\limits_{t=0}\,
 \frac{dt}{2\pi\imath\,t^{n+1}}\,
 \prod_{i=1}^{\ell+1}\prod_{m=0}^d\frac{1}{(1-tz_iq^m)}
\,+\Big(\prod_{i=1}^{\ell+1}z_i^k\Big)\,\oint\limits_{t=\infty}\,
 \frac{dt}{2\pi\imath\,t^{n+1}}\,
 \prod_{i=1}^{\ell+1}\prod_{m=0}^d\frac{1}{(1-tz_iq^m)}.
\ee   This
decomposition corresponds to the decomposition \eqref{decomptwo} \be
 \chi_G(\QM_d(\IP^\ell),\CL_k(n))\,=\,\\
=\Tr_{H^0(\QM_d(\IP^\ell),\CL_k(n))}\,\,
q^{L_0}\,e^{\sum\lambda_iH_i}+(-1)^{(\ell+1)(d+1)-1}
\Tr_{H^{(\ell+1)(d+1)-1}(\QM_d(\IP^\ell),\CL_k(n))}\,\,
 q^{L_0}\,e^{\sum\lambda_iH_i}\,.\nonumber
\ee
\end{rem}

\begin{rem}\label{remzero}
In the limit $q\to 0$ one has an integral representation for
a character $\chi_{(n,k)}^{(0)}$ of an irreducible
finite-dimensional representation
$V_{n,k,0}=\Sym^n\IC^{\ell+1}\otimes \CL_k$ of $GL_{\ell+1}$:
\be\label{qzerochar}
 \chi^{(0)}_{(n,k)}(\underline{z})\,=\,\tr_{V_{n,k,0}}
 e^{\lambda_1H_1+\ldots+\lambda_{\ell+1}H_{\ell+1}}\,=\,
 \Big(\prod_{i=1}^{\ell+1}z_i^k\Big)\,\oint\limits_{t=0}\,
 \frac{dt}{2\pi\i\, t^{n+1}}\,
 \prod_{i=1}^{\ell+1}\frac{1}{1-tz_i}.
\ee
where $z_i=\exp \lambda_i$ and 
the $GL(\ell+1)$-module $V_{n,k,0}$, $n\geq 0$ is  realized as
a zero cohomology space $H^0(\IP^\ell,\CL_k(n))$. 
\end{rem}

\section{K-theory of $\CL\IP^{\ell}_+$ and $q$-Whittaker functions}

In this section we establish a direct connection between 
$q$-deformed class one specialized $\mathfrak{gl}_{\ell+1}$-Whittaker functions 
and geometry of the space $\CL\IP^{\ell}_+$ defined as an appropriate
limit of $\CQ\CM_d(\IP^{\ell})$  when $d\to +\infty$. Geometrically
the $\CL\IP^{\ell}_+$ should be considered as a space of algebraic
disks in $\IP^{\ell}$ (see \cite{Gi1} for details). In general,   
let $LX={\rm Map}(S^1,X)$ be the space of free contractible loops in
a compact K\"ahler manifold $X$. There is a natural action
of $S^1$ on  $LX$ by loop rotations. The universal covering $\widetilde{LX}$
can be defined as a space of maps $D\to X$ of the disk $D$
considered up to a homotopy of the map preserving the image of the boundary
loop $S^1\subset D$. The group of covering transformations
of the universal cover $\widetilde{LX}\to LX$ is isomorphic to
$\pi_2(X)$. Let $\widetilde{LX_+}\subset \widetilde{LX}$ 
 be a semi-infinite cycle of loops that are
boundaries of holomorphic maps $D\to X$. For $X=\IP^{\ell}$ define 
an algebraic version $\CL \IP^{\ell}_+$ of $\widetilde{L\IP^{\ell}}_+$    
as a set of collections of regular series 
$$
 a_i(z)=a_{i,0}+a_{i,1}z+a_{i,2}z^2+\cdots,\hspace{3cm} 0\leq i\leq\ell,
$$
modulo the action of $\IC^*$. The topology on this space 
should be defined by considering $\CL \IP^{\ell}_+$ 
as a   limit of $\CQ\CM_d(\IP^\ell)$ when $d\to\infty$. This space inherits the
action of $G=S^1\times U(\ell+1)$ defined previously on 
$\CQ\CM_d(\IP^\ell)$. In the following we do not define appropriate  
topology rigorously leaving this for another occasion. Instead we
define the limit $d\to +\infty$ on the level of cohomology algebra 
$H^*(\CQ\CM_d(\IP^{\ell}),\IC)$ and the space of holomorphic sections 
of line bundles on $\CQ\CM_d(\IP^{\ell})$. 
Let us take the limit $d\to +\infty$ of
the character $A^{(d)}_{n,k}(\underline{z},q)$ given by the integral
expression \eqref{finited}. The limit of
$A^{(d)}_{n,k}(\underline{z},q)$ can be interpreted as a character of a 
$\IC^*\times GL_{\ell+1}$-module $\CV_{n,k,\infty}$ defined as
follows. Let $\CV_{k,\infty}$ be a linear space of polynomials of infinite number of
variables  $a_{i,m}$, $i=0,\ldots ,\ell$, $m\in \IZ_{\geq 0}$. 
Let $L_0$ be a generator of $\Lie(\IC^*)$, 
$T\in GL_{\ell+1}$ be a Cartan torus and 
$H_1,\ldots,H_{\ell+1}$ be a basis in $\Lie(T)$.
Define the action of $L_0$ and $H_j$ on the generators $a_{i,m}$ as
follows  
$$
L_0:\,\,a_{i,m}\longrightarrow m\,\,a_{i,m};
$$
$$
e^{\sum_j\lambda_j H_j}:\,\,a_{i,m}\longrightarrow
e^{\lambda_i}\,a_{i,m}
$$  
This supplies $\CV_{k,\infty}$ with the structure of a $\IC^*\times
GL_{\ell+1}$-module. Now 
the linear subspace $\CV_{n,k,\infty}\subset \CV_{k,\infty}$ is
defined as a subspace 
of polynomials of the variables $a_{i,m}$, $i=0,\ldots ,\ell$, $m\in \IZ_{\geq 0}$
of the total degree $n$. 

\begin{te}\label{mainTh21} Let
$\Psi^{\mathfrak{gl}_{\ell+1}}_{\underline{z}}(n,k)$
be a specialization \eqref{spec}
of the solution of $q$-deformed Toda lattice defined in the Theorem
\ref{mainth1}. Then the following holds.
\be\label{identification}
\Psi^{\mathfrak{gl}_{\ell+1}}_{\underline{z}}(n,k)\,=\,
\chi_{n,k}(\underline{z}),
\ee
where 
$$
\chi_{n,k}(\underline{z})=\lim_{d\to
\infty}\,A^{(d)}_{n,k}(\underline{z},q)=\Tr_{\CV_{n,k,\infty}}\,q^{L_0}\,e^{\sum\lambda_iH_i},
\qquad z_i=e^{\lambda_i}.
$$
\end{te}

\noindent {\it Proof}:
For the function $\chi_{n,k}(\underline{z})=\lim_{d\to
\infty}\,A^{(d)}_{n,k}(\underline{z},q)$  the following integral
representation holds: \be\label{qchar}
\chi_{n,k}(\underline{z})\,=\,
\Big(\prod_{i=1}^{\ell+1}z_i^k\Big)\,\oint_{t=0}\,
\frac{dt}{2\pi\imath\,t^{n+1}}\,
\prod_{m=1}^{\ell+1}\prod_{j=0}^\infty\frac{1}{(1-tq^jz_m)}=\ee
$$=\Big(\prod_{i=1}^{\ell+1}z_i^k\Big)\,\oint_{t=0}\,
\frac{dt}{2\pi\imath\,t^{n+1}}\, \prod_{m=1}^{\ell+1}\Gamma_q(tz_m).
$$
The relations \eqref{identification} follows directly from
the explicit integral expression \eqref{qchar} and Lemma \ref{lemma}.
The representation in terms of the trace over $\CV_{n,k,\infty}$ follows 
from the statement of Proposition \ref{refmainth} with obvious
modifications for $d\to +\infty$\, $\Box$

For $n\geq 0$ we can identify the character $A^{(d)}_{n,k}$ 
with the equivariant Euler characteristic expressed 
through  Riemann-Roch-Hirzebruch formula 
\be \chi_G({\cal{QM}}_d(\IP^\ell),{\cal{L}}_k(n))\,=\,
\int_{{\cal{QM}}_d(\IP^\ell)}
 {\rm Ch}_G({\cal{L}}_k(n))\,{\rm
 Td}_G(\CT{\cal{QM}}_d(\IP^\ell)).
\ee 
Taking the limit $d\to +\infty$ we obtain formal  
Riemann-Roch-Hirzebruch formula for  
$\chi_G(\CL \IP^{\ell}_+,{\cal{L}}_k(n))$. 
 Using the description \eqref{Ktheory}, \eqref{pairingK}
of the equivariant K-groups of projective spaces 
and taking the limit $d\to +\infty$ in the integral representation 
of the Euler characteristic \eqref{MainChar}
one obtains the following integral representation for the equivariant
Euler  characteristic of line bundles on  $\CL\IP^{\ell}_+$ 
\be\label{RRHsemiinfinite}
\chi_G(\CL \IP^{\ell}_+,{\cal{L}}_k(n))
\,=\,-\Big(\prod_{j=1}^{\ell+1}z_j^k\Big)\oint\limits_{C}\,\frac{dt}{2\pi\imath\,
t^{n+1}}\,
\prod_{i=1}^{\ell+1}\prod_{j=0}^{\infty}\frac{1}{(1-tq^jz_i)}\,=\,\\
=\,- \Big(\prod_{j=1}^{\ell+1}z_j^k\Big)\oint\limits_{C}\frac{dt}{2\pi\i\,t^{n+1}}
 \prod_{i=1}^{\ell+1}\Gamma_q(tz_i),
\ee
 where the integration contour $C$ encircles all poles except $t=0$ and   
$$
\Gamma_q(y)\,=\,\prod_{n=0}^\infty\frac{1}{1-yq^n}.
$$

\begin{prob}\label{conjNum}
Define an equivariant (co)homology theory for $\CL\IP^{\ell}_+$ in
such a way that Chern and Todd classes ${\rm Ch}_G(\CL_k(n))$, 
${\rm Td}_G(\CT\CL\IP^{\ell}_+)$ make sense and the expression   
$$
\int_{\CL\IP^{\ell}_+}
 {\rm Ch}_G({\cal{L}}_k(n))\,{\rm Td}_G(\CT\CL\IP^{\ell}_+)
$$
is  well-defined and is equal to  
\be\label{conjKtheory}
\Psi^{\mathfrak{gl}_{\ell+1}}_{\underline{z}}(n,k)\,=\,-
 \Big(\prod_{j=1}^{\ell+1}z_j^k\Big)\oint\limits_{C}\frac{dt}{2\pi\i\,t^{n+1}}
 \prod_{i=1}^{\ell+1}\Gamma_q(tz_i).
\ee
\end{prob}

\begin{rem} The conjectural relation  above provides a  description
  of the specialized $q$-deformed $\mathfrak{gl}_{\ell+1}$-Whittaker
  function as a semi-infinite period 
 \be
\Psi^{\mathfrak{gl}_{\ell+1}}_{\underline{z}}(n,k)=
\int_{\CL\IP^{\ell}_+}
 {\rm Ch}_G({\cal{L}}_k(n))\,{\rm Td}_G(\CT\CL\IP^{\ell}_+),\qquad
 n\geq 0.
\ee
\end{rem}
The $K$-theory of the semi-infinite  spaces $\CL\IP^{\ell}_+$ 
is closely connected with a quantum version of $K$-theory of projective
 spaces  proposed in \cite{GiL}.  
The generating function $F(n,\underline{z},q)$ of the correlation
functions in $K$-theory version of Gromov-Witten theory with the target
space $\IP^{\ell}$  obeys the following difference equation 
\cite{GiL}
\be\label{ktheory}
 \Big\{\prod_{i=1}^{\ell+1}(1-z_iT^{-1})\Big\}\cdot
 F(n,\underline{z},q)\,=\,q^nF(n,\underline{z},q),
\ee
where $T\cdot f(n)=f(n+1)$.
The specialized $q$-deformed $\mathfrak{gl}_{\ell+1}$-Whittaker
function satisfies  the same equation  \eqref{ktheory}
( see Lemma \ref{lemma} and relation
\eqref{eqdif}). Therefore the
Whittaker function  can be considered as a correlation function 
of some special operator singled out by the class one condition
(i.e. the condition $\Psi^{\mathfrak{gl}_{\ell+1}}_{\underline{z}}(n,k)=0$ for $n<0$).
We provide some information on this operator in the Section \ref{S1loc}.

\section{Quantum  cohomology and Whittaker function}

In the previous Section we proposed 
a description of  $q$-deformed  class one specialized 
$\mathfrak{gl}_{\ell+1}$-Whittaker function  
in terms of a semi-infinite version of Riemann-Roch-Hirzebruch 
theorem. This expresses the $q$-Whittaker function as  a semi-infinite
period. Its  classical (i.e. non-deformed) counterpart can be also 
expressed in terms of a semi-infinite period. 
In this Section we provide this conjectural representation
\footnote{Whittaker
  functions naturally arise in the description of Gromov-Witten  
invariants of flag spaces. In the mirror dual description 
they expressed in terms of periods of  top-dimensional
holomorphic forms on non-compact Calabi-Yau spaces \cite{Gi2}. 
Thus, the possibility to express Whittaker functions as semi-infinite periods  
leads to a formulation  of the  mirror symmetry as an identification of two
period maps - semi-infinite and finite ones.}.  

We start from recalling the notion of quantum cohomology. 
The  quantum cohomology $QH^*(X)$ of a compact symplectic manifold $X$
 can be defined in terms of
semi-infinite geometry of a universal cover $\widetilde{LX}$ of the loop
space $LX$. One of the descriptions 
is given by a Morse-Smale-Bott-Novikov-Floer complex
constructed  in terms of
critical points of an area functional on $\widetilde{LX}$. Its
cohomology groups (interpreted as Floer cohomology groups 
$FH^*(\widetilde{LX})$  of $\widetilde{LX}$) are  isomorphic to
 the semi-infinite cohomology $H^{\infty/2+*}(LX)$ arising naturally
in the Hamiltonian formalism of a topological two-dimensional sigma
model with the target space $X$.   In the following we will use an
equivariant version of quantum cohomology $QH^*(\IP^{\ell})$
of projective spaces considered
in \cite{Gi1} (see also \cite{CJS} for a non-equivariant version).

We have defined the universal covering $\widetilde{LX}$
of the loop space $LX$ as a space of maps $D\to X$ of the disk $D$
considered up to a homotopy preserving the image of the boundary
loop $S^1\subset D$. The group of covering transformations
of the universal cover $\widetilde{LX}\to LX$ is isomorphic to
the image $\Gamma\subset H_2(X)$ of the Hurewicz homomorphism
$\pi_2(X)\to H_2(X)$ where $H_2(X)$ denotes integral homology modulo
torsion. Let $r$ be the rank of $\Gamma$ and
$\IC\bigl[\Gamma \bigr]\simeq\IC[u_1^{\pm1},\ldots,u_r^{\pm1}]$
be its group algebra.

As  a vector space  the quantum cohomology $QH^*(X)$  
of $X$ as a vector space is  isomorphic
to the ordinary cohomology $H^*(X,\,\IC[u_1^{\pm1},\ldots,u_r^{\pm1}])$,
over the group algebra
$\IC\bigl[\Gamma \bigr]\simeq\IC[u_1^{\pm1},\ldots,u_r^{\pm1}]$.
Let $S^1$ act on the loop space $LX$ by loop rotations.
For the corresponding $S^1$-equivariant quantum cohomology
we have the following isomorphism of vector spaces: 
$$
QH^*_{S^1}(X)\,=\,
 H^*\bigl(X,\IC[u_1^{\pm1},\ldots,u_r^{\pm1}](\hbar)\,\,\bigr),
$$
where we use the identification
$$
H^*_{S^1}(\pt,\IC)=H^*(BS^1,\IC)=\IC[\hbar],
$$
and the standard localization of the equivariant cohomology 
$H^*_{S^1}(\pt)$ with respect to the maximal ideal generated by $\hbar$ is implied.

The quantum cohomology space $QH^*_{S^1}(X)$ has a natural structure of
a module over an algebra $\CD$ generated by $u_i=\exp\tau_i$,
$v_i=-\hbar\pr/\pr\tau_i$, $i=1,\ldots,r$.  
More precisely, $QH^*_{S^1}(X)$ as a linear space over $\IC(\hbar)$ 
is generated by solutions of the system of linear differential equations 
\be\label{flatcon}
\nabla_i\,f(\underline{\tau})=0 ,\qquad
f(\underline{\tau})=(f_1(\underline{\tau}),f_2(\underline{\tau}),
\ldots,f_n(\underline{\tau})),\qquad n=\dim H^*(X),
\ee
where the  flat connection $\nabla=\sum_{i=1}^rd\tau_i\nabla_i$
provides an action of $v_i$ on $QH^*_{S^1}(X)$.

The $\CD$-module $QH^*_{S^1}(\IP^{\ell})$ is of rank one. It is generated 
over $\CD$ by an element $f_*$ satisfying the relation
$(v^{\ell+1}-u)f_*=0$ i.e. the  quantum cohomology can be
represented as $QH^*(\IP^{\ell})\simeq \CD/(v^{\ell+1}-u)$.
Explicitly we have the differential equation for the generator $f_*(\tau)$      
\be\label{givpl}
\Big\{\,\Big(-\hbar\frac{\pr}{\pr\tau}\Big)^{\ell+1}\,-\,
e^{\tau}\,\Big\}\,f_*(\tau,\hbar)\,\,=\,\,0.\ee
The representation  \eqref{flatcon} arises after 
transformation of \eqref{givpl} to the matrix 
differential equation of the first order. 

The $(S^1\times U_{\ell+1})$-equivariant analog 
 of quantum cohomology $QH^*(\IP^{\ell})$ 
allows for a similar representation with  the differential
equation  \eqref{givpl}  replaced by
 \be\label{poneeq}
\Big\{\,\prod_{k=1}^{\ell+1}\Big(-\hbar\frac{\pr}{\pr
\tau}-\lambda_k\Big)\,-\,e^{\tau}\,\Big\}\,\,f_*(\tau,
\underline{\lambda},\hbar)\,\,=\,\,0.\ee

\begin{lem}
The general solution of \eqref{poneeq} is given by
a linear combination of the integrals
\be\label{integthree}
 f^{(a)}(\tau,\underline{\lambda},\hbar)=\int_{\gamma_a} d\lambda\,\,
 e^{-\frac{\lambda\tau}{\hbar}} \,\,
 \prod_{k=1}^{\ell+1}\,\hbar^{\frac{\lambda-\lambda_k}{\hbar}}\,
 \Gamma\Big(\frac{\lambda-\lambda_k}{\hbar}\Big),
\qquad a=1,\cdots ,n
\ee
with a suitable choice of integration contours $\gamma_a$.
\end{lem}
\noindent {\it Proof}:
Note  that the function
\be
 Q(\lambda,\underline{\lambda})\,=\,
 \prod_{k=1}^{\ell+1}\,\hbar^{\frac{\lambda-\lambda_k}{\hbar}}\,
 \Gamma\Big(\frac{\lambda-\lambda_k}{\hbar}\Big).
\ee
obeys the difference equation
\be
\prod_{k=1}^{\ell+1}\bigl(\lambda-\lambda_k\bigr)\,
Q(\lambda,\underline{\lambda})\,=\,
Q(\lambda+\hbar,\underline{\lambda}).
 \ee
Therefore, the function
 \be\label{classonesol} 
 f(\tau,\underline{\lambda})=\int_{\gamma} d\lambda\,\,
 e^{-\frac{\lambda\tau}{\hbar}} \,\,
 Q(\lambda,\underline{\lambda},\hbar)\ee
satisfies  \eqref{poneeq} provided the choice of the contour $\gamma$ 
allows for an integration by parts. 
The contours can be  chosen is such a way that the total
derivatives do not give a contribution into the integral
\eqref{integthree} 
$\Box$

A particular choice of $\gamma$ in \eqref{classonesol} gives us a
special solution of the equation \eqref{poneeq} 
\be\label{f*}
 f_*(\tau,\underline{\lambda},\hbar)\,=\,
\int\limits_{\sigma-\imath \infty}^{\sigma+\imath\infty}\,
d\lambda\,\,\,e^{-\frac{\lambda\tau}{\hbar}}\,\,\,\prod_{k=1}^{\ell+1}
\hbar^{\frac{\lambda-\lambda_k}{\hbar}}\, \Gamma
\Big(\frac{\lambda-\lambda_k}{\hbar}\Big), 
\ee
where $\sigma$ is such that $\sigma< min\,\{\,{\rm
  Re}\,\lambda_j,j=1,\ldots,\ell+1\}$. 
This is a unique solution of \eqref{poneeq} exponentially decaying 
when $\tau\to +\infty$. 

\begin{rem}
In the case of  $\ell=1$,  the differential equation \eqref{givpl} 
is equivalent to an eigenvalue problem for the quadratic Hamiltonian of the 
$\mathfrak{sl}_2$-Toda chain. 
The solution given in the integral form  \eqref{f*} coincides in this
case with the Mellin-Barnes
representation of the $\mathfrak{gl}_2$-Whittaker function
(\ref{mellin2}) for $x_2=0$ and $x_1=\tau$.
\end{rem}

Replacing formally $\Gamma$-functions by infinite products 
over their poles one has  for $ f_*(\tau,\underline{\lambda},\hbar)$ 
the following expression 
 \be\label{formsolH}
 \int\,dx\,e^{\tau x/\hbar}\,\,
 \prod_{j=1}^{\ell+1}\prod_{n=0}^\infty\frac{1}{x-\lambda_j-\hbar
 n}.
\ee 
This formal representation can be interpreted using the model for
cohomology of $\CQ\CM_d(\IP^{\ell})$ discussed n Section
\ref{modelcoh}. Naively \eqref{formsolH} can be considered as an
integral over $\CL\IP^{\ell}_+$ of $\exp(\tau\omega/\hbar)$
where $\omega$ is an element of the second  $S^1\times
U_{\ell+1}$-equivariant cohomology of $\CL\IP^{\ell}_+$. Recall that
we define  $\CL\IP^{\ell}_+$ 
on the level of cohomology as a limit of $\CQ\CM_d(\IP^{\ell})$ when 
$d\to +\infty$. However a correct regularization for 
 \eqref{formsolH} is given by \eqref{f*} and, thus, a geometric
 interpretation of \eqref{f*} implies some modification of
 $\CL\IP^{\ell}_+$. We attribute the difference
 between \eqref{formsolH} and \eqref{f*} to the fact that 
the proper interpretation of $\CL\IP^{\ell}_+$ as $d\to +\infty$ limit
deserves more care in this case and does not coincide with
a straightforward limit $d\to +\infty$ on the level of cohomology. 
Let us denote the corresponding hypothetically modified limit by ${\bf L}\IP^{\ell}_+$.
 
\begin{prob} Find the space ${\bf L}\IP^{\ell}_+$ and construct
  equivariant (co)homology for ${\bf L}\IP^{\ell}_+$ in such a way
  that the integral  
$$ 
\int_{{\bf
    L}\IP^{\ell}_+}\,e^{\tau\,\omega/\hbar},
\qquad \omega\in H_{S^1\times U_{\ell+1}}^2({\bf L}\IP^{\ell}_+,\IC)
$$
is well-defined and is equal to $f_*(\tau,\underline{\lambda},\hbar)$ 
given by \eqref{f*}. \end{prob}

\section{$S^1$-localization }\label{S1loc}

In this Section we calculate  the equivariant  
Euler characteristic  $\chi_{G}(\QM_d(\IP^{\ell}),\CL_k(n))$
for $G=S^1\times U_{\ell+1}$ using 
Borel localization for $S^1$-action. This yields a direct
relation between our construction of $q$-Whittaker functions 
and the results of \cite{GiL}.

The character \eqref{MainChar} can be  calculated using an
equivariant localization as follows. 
We have  a compact Lie group $G=S^1\times U_{\ell+1}$ acting  
on a  projective space $X=\CQ\CM_d(\IP^{\ell})\equiv
\IP^{(\ell+1)(d+1)-1}$. Recall that $\CQ\CM_d(\IP^{\ell})$ is defined
as a set of $(\ell+1)$ polynomials 
each of degree $d$ considered up to common constant factor
 \be\nonumber
\bigl(a_0(y),a_1(y),\ldots
a_\ell(y)\bigr)\,\sim \,\bigl
(\rho\,a_0(y),\rho\,a_1(y),\ldots
\rho\,a_\ell(y)\bigr),\qquad \rho\in \IC^*,
\ee
where
$$
a_k(y)=\sum_{i=0}^d\,a_{k,i}\,y_1^i\,y_2^{d-i}, \qquad k=0,\ldots,\ell.
$$
The action of an element $q=e^{\imath
  \theta}$ of $S^1$ on $a_{k,j}$ is given by 
$$
q:\,\, a_{k,j}\rightarrow q^j a_{k,j}.
$$
The line bundles $\CL_k(n)$ are equivariant with
respect to the action of $S^1$.
Let  $X^{S^1}\subset X$ be a set of $S^1$-fixed points. It is a union of
smooth components $Y_i$. The Bott localization formula
gives the following expression for the equivariant Euler characteristic
\eqref{MainChar} (see e.g. \cite{BGV}) 
\be\label{IndexTh}
 \chi_{G}(\QM_d(\IP^{\ell}),\CL_k(n))=\sum_{Y_i\in X^{S^1}}
\,\int_{Y_i}\frac{{\rm Ch}_G(\CL_k(n)|_{Y_i}){\rm Td}_G(TY_i)}{{\rm E}_G(\CN_{Y_i}) },
\ee
where the sum  runs over all components $Y_i$ in $X^{S^1}$,
$\CN_{Y_i}$ is the  normal bundle of $Y_i$ in $X$,
${\rm Ch}_G(\CL_{k}(n))$, ${\rm Td}_G(Y_i)$  are  equivariant
Chern character and Todd class, and ${\rm E}_G(\CN_{Y_i})$ is the equivariant
Euler class of $\CN_{Y_i}$.

It is easy to infer that the subvarieties $Y_i$, $i=0,\ldots
,\ell$  in $\CQ\CM_d(\IP^{\ell})$  are isomorphic to the projective spaces
$\IP^{\ell}$ and are defined by the equations
\be\nonumber
Y_i=\{a_{k,j}=0,\, j\neq i\},\qquad i=0,\ldots ,\ell.
\ee
To calculate the action of $S^1$ on the normal bundle to  $Y_{i}$ we
consider the intersection of $Y_i$ with
open subsets $U_{a_{k,i}}=\{a_{k,i}\neq 0\}$ of $\CQ\CM_d(\IP^{\ell})$.
Natural coordinates on $U_{a_{k,i}}$ are
$$
\xi_{r,j}=a_{r,j}/a_{k,i},\qquad (r,j)\neq (k,i),
$$
and the intersections $Y_i\cap U_{a_{k,i}}=\{a_{k,i}\neq 0\}$
are defined by the equations
$$
\xi_{r,i}=0,\qquad r\neq k.
$$
Thus one can take a collection of coordinates $\xi_{r,j}$, $j\neq i$
as a local section of the dual to the normal bundle $\CN_{Y_i}$.
The action of $S^1$ on $\CN_{Y_i}$ can  then be found by considering
the action on section $\xi_{r,j}$:
$$
\xi_{r,j}\rightarrow q^{j-i}\xi_{r,j}.
$$
Similarly,  one can show that
$q\in S^1$ acts on the restriction of the line bundle $\CO()$  
on $Y_i$ by  multiplication  on $q^{n\,i}$. The fixed point formula 
\eqref{IndexTh} reduces to the following explicit expression 
\be\label{decom}
\chi_{G}(\QM_d(\IP^{\ell}),\CL_k(n))=\\
=\,-\,\Big(\prod_{j=1}^{\ell+1}\,z_j^k\Big)\, \sum_{i=0}^{d}
\int_{C_0}\,\frac{dt}{2\pi\imath t^{n+1}}\,
\frac{q^{ni}}{\prod_{j=1}^{\ell+1}\prod_{m=0,m\neq
i}^{d}(1-tz_jq^{m-i})} \frac{1}{\prod_{j=1}^{\ell+1}(1-tz_j)}\nonumber,
\ee
where the integration contour $C_0$ encircles the $(\ell+1)$  poles defined by the
equations $t=z_j^{-1}$, $j=1,\ldots ,\ell+1$.

\begin{lem}
The following identity holds for $n\geq 0$
\be
 \oint_{C}\,\frac{dt}{2\pi\imath t^{n+1}}\,
\frac{1}{\prod_{j=1}^{\ell+1}\prod_{m=0}^{d}(1-tz_jq^m)}=\nonumber
\ee
\be
=\sum_{i=0}^{d}
\int_{C_0}\,\frac{dt}{2\pi\imath t^{n+1}}\,
\frac{q^{ni}}{\prod_{j=1}^{\ell+1}\prod_{m=0,m\neq
i}^{d}(1-tz_jq^{m-i})} \frac{1}{\prod_{j=1}^{\ell+1}(1-tz_j)},
\ee
where the integration contour $C$ encircles  poles defined
by the equations $t=z_j^{-1}q^{-i}$, $j=1,\ldots,\ell+1$, $i=0,\ldots ,d$
and the integration contour $C_0$ encircles $(\ell+1)$ poles defined
by the equations $t=z_j^{-1}$, $j=1,\ldots,\ell+1$. 
\end{lem}
\noindent{\it Proof}:
We have that 
$$
\oint_{C}\,\frac{dt}{2\pi\imath t^{n+1}}\,
\frac{1}{\prod_{j=1}^{\ell+1}\prod_{m=0}^{d}(1-tz_jq^m)}
\,=\,\sum_{i=0}^{d}
\int_{C_i}\,\frac{dt}{2\pi\imath t^{n+1}}\,
\frac{1}{\prod_{j=1}^{\ell+1}\prod_{m=0}^{d}(1-tz_jq^m)},
$$
where the integration contour $C_i$ encircles $(\ell+1)$ poles defined
by the equations $t=z_j^{-1}q^{-i}$, $j=1,\ldots,\ell+1$. Making the change of
variables $t\to tq^{-i}$ in the r.h.s.,  we obtain that 
$$
\oint_{C}\,\frac{dt}{2\pi\imath t^{n+1}}\,
\frac{1}{\prod_{j=1}^{\ell+1}\prod_{m=0}^{d}(1-tz_jq^m)}=
$$
$$
\sum_{i=0}^{d} \int_{C_i}\,\frac{dt}{2\pi\imath t^{n+1}}\,
\frac{1}{\prod_{j=1}^{\ell+1}\prod_{m=0}^{d}(1-tz_jq^m)}=\sum_{i=0}^{d}q^{ni}
\int_{C_0}\,\frac{dt}{2\pi\imath t^{n+1}}\,
\frac{1}{\prod_{j=1}^{\ell+1}\prod_{m=0}^{d}(1-tz_jq^{m-i})}
$$
$\Box$

We are going to consider  a continuation of the expression
\eqref{decom}  to $q\in \IC^*$, $|q|<1$ and the limit $d\to \infty$.

\begin{prop} \label{GiLrep}
The specialization \eqref{spec}, \eqref{qWhit}
of the $q$-deformed $\mathfrak{gl}_{\ell+1}$-Whittaker
function can be written in the following form 
\be
\Psi^{\mathfrak{gl}_{\ell+1}}_{\underline{z}}(n,k)\, =
\<I_{n,k}(\underline{z})\,\tilde{L}(\underline{z}),[\IP^{\ell}]\>_K,
\ee
where
$$
\tilde{L}(\underline{z},t)=
\frac{1}{\prod_{j=1}^{\ell+1}\prod_{k=1}^{\infty}(1-tz_jq^{k})}=
\prod_{j=1}^{\ell+1}\Gamma_q(qtz_j),
$$
\be\label{GLigen}
I_{n,k}(\underline{z},t)=\,\Big(\prod_{j=1}^{\ell+1}\,z_j^k\Big)\,
t^{-n}\,\sum_{i=0}^{\infty}q^{ni}
\frac{1}{\prod_{j=1}^{\ell+1}\prod_{m=1}^{i}(1-tz_jq^{-m})},
\ee
and the pairing $\<\,,\,\>_K$ is the standard pairing
\eqref{pairingK} on $K_{U_{\ell+1}}(\IP^{\ell})$ taking
values in $K_{U_{\ell+1}}(pt)$.
\end{prop}

The representation of the  Whittaker function given in Proposition \ref{GiLrep}
establishes a direct connection with the results of  Givental-Lie \cite{GiL}.
In \cite{GiL} the function \eqref{GLigen} was interpreted as a universal solution of
the reduction \eqref{eqdif} of $q$-deformed
$\mathfrak{gl}_{\ell+1}$-Toda chain. 
Indeed, the function $I_{n,k}(\underline{z},t)$  satisfies the 
eigenvalue problem
\be\label{eigprop}
\prod_{j=1}^{\ell+1}\,(1-z_iT^{-1})\,
I_{n,k}(\underline{z},t)\,=\,q^n\,
I_{n,k}(\underline{z},t),
\ee modulo the relation
$\prod_{j=1}^{\ell+1}(1-tz_j)=0$ holding  in $K_{U_{\ell+1}}(\IP^{\ell})$ and  is  uniquely
determined by the  normalization condition
$$
I_{n,k}(\underline{z},t)|_{q=0}\,=\,\Big(\prod_{j=1}^{\ell+1}\,z_j^k\Big)\,\,t^{-n},
\qquad n\geq 0. 
$$
The solution $I_{n,k}(\underline{z},t)$ is universal in the sense that
taking the pairing 
\be\label{welldefpair}
\<I_{n,k}(\underline{z}),f\>_K\,\,=\,\,
 -\frac{1}{2\pi\i}\oint\limits_{C_0}\!\! \frac{dt}{t}\,\,
 \frac{I_{n,k}(\underline{z},t)\,f(t)}{\prod_{j=1}^{\ell+1}(1-z_jt)},
\ee
with arbitrary $f\in K_{U_{\ell+1}}(\IP^\ell)$ 
one obtains a solution \eqref{welldefpair} 
of the $q$-deformed reduced 
$\mathfrak{gl}_{\ell+1}$-Toda chain \eqref{eqdif}.

\section{Semi-infinite Todd genus and $q$-Gamma function}

In the explicit expression  
for the cohomological pairing on $\CL\IP^{\ell}_+$ conjectured in Problem \ref{conjNum}  
the $q$-Gamma function $\Gamma_q$ plays the role similar 
to the Todd genus in the analogous pairing for underlying finite-dimensional
space $\IP^{\ell}$. The $S^1$-localization 
discussed in the previous section reduces
the pairing of the Chern  and Todd classes on $\CL\IP^{\ell}_+$  
to the paring of some cohomology classes on $\IP^{\ell}$. 
It is an interesting problem to interpret the resulting cohomology
classes on $\IP^{\ell}$ in terms of some geometric objects on $\IP^{\ell}$. 
For example in an analogous case of
$S^1$-localization of $K$-theory on the loop space $LX$, the elliptic 
genus of $X$ arises. The corresponding elliptic cohomology is an
instance of an extraordinary cohomology theory. 
In this Section we discuss the result of $S^1$-localization on 
$\CL\IP^{\ell}_+$ from the  extraordinary cohomology perspective. 
We demonstrate that  an intrinsic 
non-local nature of the $\IC^*$-localization on $\CL\IP^{\ell}_+$ 
obstructs a  straightforward  relation with the formalism 
of extraordinary cohomology theories and corresponding multiplicative 
genera. Let us remark that  classical  $\Gamma$-function 
was considered as a candidate for a topological genus by Kontsevich in \cite{K}.
Also a kind of $\Gamma$-genus also appeared in obviously related
context in \cite{Li}, \cite{Ho} (see also \cite{CGi1}, \cite{CGi2} 
for a discussion of formal groups in a  quantum version of cobordism theory).

We first recall  standard facts on multiplicative
topological genera and formal group laws corresponding to 
complex oriented cohomology theories (see e.g. \cite{BMN}).
The Hirzerbruch multiplicative genus is a homomorphism
$\varphi:\,\Omega^*\to \CR$ of the  ring of complex cobordisms
$\Omega^*=\Omega^*({\rm pt})$ to  a ring of coefficients $\CR$.
One has a Thom isomorphisms $\Omega^*\otimes
\IQ=\IQ[x_1,x_2,\cdots]$, $\deg(x_i)=-2i$ and the topological genus
$\varphi$ is characterized by its values on generators $x_i$ that can be
represented by Pontryagin-Thom duals to complex projective spaces $\IP^{i}$. Equivalently
$\varphi$  defined over $\IQ$ can be  described  in terms of
 a one-dimensional commutative formal group law
\be\label{fglaw}
f_{\varphi}(z,w)=e_{\varphi}(\log_{\varphi}(z)+\log_{\varphi}(w)),
\ee
expressed through the  logarithm function $$
\log_{\varphi}(z)=z+\sum_{n=1}^{\infty}\frac{\varphi([\IP^{n}])}{n+1}\,z^{n+1},
$$
and its inverse $e_{\varphi}(u)$.
For instance rational cohomology and $K$-theory correspond to additive and
multiplicative group laws
$$
f_H(z,w)=z+w,\qquad f_K(z,w)=z+w-zw.
$$
To a genus $\varphi$ one associates a multiplicative sequence $\{\Phi_n(c_i)\}$,
$\deg(\Phi_n)=n$ of cohomology classes  
$$
P_{\varphi}=\sum_{n=0}^{\infty}
\Phi_n(c_i)=\prod_{j=1}^N\frac{x_j}{e_{\varphi}(x_j)},
$$
and a map
$$
X\rightarrow \varphi(X)=\<\Phi(\CT X),[X]\>,\hspace{1.5cm}
\dim_{\IC}\,X=N.
$$
Here $\CT X$ is the tangent bundle to a manifold $X$, $[X]$ is the 
fundamental class in the homology of $X$ and $\<\,,\,\>$ is a standard
pairing.  The classes $x_i$ are defined in terms of Chern classes
$c_i$ of $\CT X$  using a splitting of $\CT X$
$$
c(X)=1+\sum_{i=1}^N c_i(X)=\prod^N_{j=1}\,(1+x_j).
$$
In the case of additive and multiplicative group laws  we have
respectively that 
$$
 P_{\varphi_H}(x)=1,\qquad  \log_{\varphi_H}(z)=z,\qquad
 e_{\varphi}(u)=u;
$$

$$
P_{\varphi_K}(x)\,=\,\prod_{j=1}^n\frac{x_j}{1-e^{-x_j}}, \qquad
\log_{\varphi_K}(z)=-\ln(1-z), \qquad
e_{\varphi_K}(u)\,=\,1-e^{-u}\,.
$$
Note that $P_{\varphi_K}(x)$ defines the Todd class of $\CT X$. For
example, the 
equivariant Riemann-Roch-Hirzebruch theorem for a line bundle $\CL_k(n)$ on
$\IP^{\ell}$ can be represented in the following form 
 \bqa\label{Ktheoryzero}
 \chi_{U_{\ell+1}}(\IP^{\ell},\CL_{k}(n))\,=\,
  \bigl\langle{\rm Ch}_{U(\ell+1)}\bigl(\CL_k(n)\bigr)\,
  {\rm Td}_{U(\ell+1)}(\CT\IP^\ell),\,[\IP^\ell]
  \bigr\rangle\\ \nonumber\,=\,
  \frac{1}{2\pi\imath}\oint\limits_{C_0}\!dx\,\,
  e^{nx+k(\lambda_1+\ldots+\lambda_{\ell+1})}\,
  \prod_{i=1}^{\ell+1}\frac{1}{1-e^{\lambda_i-x}}\\ \nonumber
  \,=\, \frac{1}{2\pi\imath}\oint\limits_{C_0}dx\,\,
  e^{nx+k(\lambda_1+\ldots+\lambda_{\ell+1})}\,
  \prod_{i=1}^{\ell+1}\frac{1}{e_{\varphi_K}(x-\lambda_i)}. 
 \eqa
Here $e_{\varphi_K}(x)$ is the exponent corresponding to $K$-theory 
(see the Remark \ref{remzero}).

Using the conjectural relation in Problem \ref{conjNum}, 
the $S^1\times U_{\ell+1}$-equivariant Riemann-Roch-Hirzebruch theorem
for a trivial line bundle
on $\CL\IP^{\ell}_+$ can be represented in the form
similar to  \eqref{Ktheoryzero} for $k=n=0$
\be\label{SiKtheory}
 \chi_{S^1\times U_{\ell+1}}(\CL\IP^{\ell}_+,\CL_{0}(0))\,=\,
 \frac{1}{2\pi\imath}\oint\limits_{C}\!dx\,\,
 \prod_{j=1}^{\ell+1}\prod_{m=0}^\infty
 \frac{1}{1-e^{\lambda_j+m\hbar-x}}=
\ee
$$
=\frac{1}{2\pi\imath}\oint\limits_{C}\!dx\,\,
\prod_{i=1}^{\ell+1}\frac{1}{e_{\varphi_q}(x-\lambda_i)},
$$
where $C$ encircles all the poles $x=\lambda_j+m\hbar$,
$j=1,\ldots,(\ell+1)$, $m\in \IZ_{\geq 0}$ and  we used a notation 
\be\label{Kgen}
 e_{\varphi_q}(u;\hbar)=\frac{1}{\Gamma_q(e^{-u})}\,=\,
 \prod_{n=0}^\infty\,\bigl(1-e^{n\hbar-u}\bigr)\,.
\ee
However, despite a similarity of \eqref{Ktheoryzero}
and \eqref{SiKtheory} the difference in integration contours does not
allow directly to interpret $e_{\varphi_q}(u;\hbar)$ as a topological genus corresponding 
to a extraordinary cohomology theory on $\IP^{\ell}$. The way to transform 
\eqref{SiKtheory} into an integral over the contour $C_0$ was 
discussed in Section \ref{S1loc}. As a result the 
integration contour in \eqref{SiKtheory} can be replaces by $C_0$
at the expense of multiplying the integrand by the 
correction factor $I_{0,0}(\underline{z},e^{-x})$ (see \eqref{GLigen}
for explicit expression for  $I_{n,k}(\underline{z},e^{-x})$). Now we have an expression
for the equivariant Euler characteristic on $\CL\IP^{\ell}_+$ in terms
of the pairing of cohomology classes on $\IP^{\ell}$. 
However this correction factor appears to spoil the multiplicative property of 
\eqref{Kgen}. The underlying reason for this is 
the appearance of an infinite number of copies of $\IP^{\ell}$
as components of fixed point set in $S^1$-localization.
Thus, the situation is very much different from,
for example, the elliptic genus (see e.g. \cite{Se})
where the fixed point set of $S^1$ acting on $LX$ is simply $X$ itself.
It is conceivable that the failure to interpret $\S^1$-localization
on $\CL\IP^{\ell}_+$ in terms of an extraordinary topological genus
implies actually the existence of a meaningful quantum version of the 
an extraordinary cohomology theory formalism.

\vskip 1cm

\noindent {\small {\bf A.G.}: {\sl Institute for Theoretical and
Experimental Physics, 117259, Moscow,  Russia; \hspace{8 cm}\,
\hphantom{xxx}  \hspace{4 mm} School of Mathematics, Trinity
College, Dublin 2, Ireland; \hspace{8 cm}\,
\hphantom{xxx}   \hspace{3 mm} Hamilton
Mathematics Institute, TCD, Dublin 2, Ireland;}}

\noindent{\small {\bf D.L.}: {\sl
 Institute for Theoretical and Experimental Physics,
117259, Moscow, Russia};\\
\hphantom{xxxxxx} {\it E-mail address}: {\tt lebedev@itep.ru}}\\

\noindent{\small {\bf S.O.}: {\sl
 Institute for Theoretical and Experimental Physics,
117259, Moscow, Russia};\\
\hphantom{xxxxxx} {\it E-mail address}: {\tt Sergey.Oblezin@itep.ru}}
\end{document}